\documentclass[12pt]{article}

\usepackage{amstext}
\usepackage{latexsym}
\usepackage{amssymb}
\usepackage{amsmath}
\usepackage{tensor} 

\def\R{\mathbb R}

\def\N{\mathbb N}

\newcommand{\eps}{\varepsilon}

\begin{document}

\newtheorem{theorem}{Theorem}[section]
\renewcommand{\thetheorem}{\arabic{section}.\arabic{theorem}}
\newtheorem{definition}[theorem]{Definition}
\newtheorem{deflem}[theorem]{Definition and Lemma}
\newtheorem{lemma}[theorem]{Lemma}
\newtheorem{example}[theorem]{Example}
\newtheorem{remark}[theorem]{Remark}
\newtheorem{remarks}[theorem]{Remarks}
\newtheorem{cor}[theorem]{Corollary}
\newtheorem{pro}[theorem]{Proposition}
\newtheorem{proposition}[theorem]{Proposition}

\renewcommand{\theequation}{\thesection.\arabic{equation}}
 
%
% 2) check labels and references
%
% 3) dies ist die Version mit Faktor zu $\chi$ 
%
% 4) \R^3 bei \int oder nicht; \lambda_{f_n} oder \lambda_n etc.; 
%
% 5) include $\alpha$? 
%
% 7) Bedingung an $k$ 
%
% 8) alle Aussagen checken, wo irgendwas \in {\cal A}; 
%    Variante mit anderer Definition;  
%
% 11) Notation konsistent mit erster Arbeit: {\cal D} und {\cal H}, 
%     $\hat{\Phi}$ und $\chi$, $\sigma$ und $k$; 
%
%
%
%
%
%

\title{Existence of a minimizer to the particle number-Casimir functional 
for the Einstein-Vlasov system}

\author{H{\aa}kan Andr\'{e}asson\\
        Mathematical Sciences\\
        Chalmers University of Technology\\
        G\"{o}teborg University\\
        S-41296 G\"oteborg, Sweden\\
        email: hand@chalmers.se\\
        \ \\
        Markus Kunze\\
        Mathematisches Institut\\
        Universit\"at K\"oln\\
        Weyertal 86-90\\
        D-50931 K\"oln, Germany\\
        email: mkunze1@uni-koeln.de}      

\maketitle

\begin{abstract}

\noindent
In 2001 Wolansky \cite{Wol} introduced a particle number-Casimir functional 
for the Einstein-Vlasov system. Two open questions are associated with this functional. 
First, a meaningful variational problem should be formulated and the existence of a minimizer 
to this problem should be established. The second issue is to show that a minimizer, 
for some choice of the parameters, is a static solution of the Einstein-Vlasov system. 
In the present work we solve the first problem by proving the existence of a minimizer 
to the particle number-Casimir functional. 
On the technical side, 
it is a main achievement that we are able to bypass the non-compactness of minimizing 
sequences by new arguments in both $v$-space and $x$-space, 
which might have several further applications. 
We note that such compactness results for the Einstein-Vlasov system 
have been absent in the literature, whereas similar results have been known in the Newtonian case. 
We also provide arguments which give strong support that minimizers corresponding to small masses 
are static solutions of the Einstein-Vlasov system. 
Furthermore, our analysis leads us to propose a new stability criterion for static solutions: 
We conjecture that a static solution for which the Casimir-binding energy is positive is stable for mass-preserving perturbations. 

\end{abstract}

%%%%%%%%%%%%%%%%%%%%%%%%%%%%%%%%%%%%%%%%%%%%%%%%%%%%%%%%%%%%%%%%%%%%%%%%%%%%%%%%%%%%%%%%

\section{Introduction}
\setcounter{equation}{0}

The Einstein-Vlasov system is a PDE system from Mathematical General Relativity, 
where Einstein's field equations are used jointly with a kinetic description of the matter 
to model the time evolution of a collection of non-colliding objects (`particles'); 
the objects to be described in astrophysics include stars, galaxies or galaxy clusters. 
The system simplifies considerably, if some symmetry assumptions are added, 
like spherical symmetry.  
 
It is well-known that there exists a large number of spherically symmetric 
static (i.e., time-independent) solutions to the Einstein-Vlasov system \cite{Rn2,A2}. 
Whether or not a subclass of these solutions can be obtained as minimizers of a variational problem, 
as should be expected from general physics principles,  
has however remained an open question. The main reason why this is of interest is related to stability. 
In the Newtonian analogue, i.e., for the Vlasov-Poisson system, stability of a class of static solutions 
was first shown by Guo and Rein as a consequence of their result that the solutions are minimizers 
of an energy-Casimir functional \cite{GR}. Many quite general results on stability for the Vlasov-Poisson system 
have since then been obtained, see \cite{LMR,Rn2} and the references therein. 

For the Einstein-Vlasov system it was claimed by Wolansky in \cite{Wol} that a class of static solutions 
of the Einstein-Vlasov system are minimizers to a particle number-Casimir functional. However, the proof 
in \cite{Wol} contains serious flaws as observed in \cite{AK}, and hence the problem was still open. 
It naturally splits into two parts: Firstly, one needs to formulate a variational problem 
and to show that a minimizer to this variational problem does exist. The second part is to prove that a minimizer, 
for some choices of the parameters, is a static solution of the Einstein-Vlasov system. 

The main purpose of the present work concerns part one, i.e., we set up a variational problem 
and establish that a minimizer to this problem exists. To show existence of a minimizer to a variational problem 
often requires that compactness properties of the minimizing sequence can be derived. Indeed, this is central 
in the Vlasov-Poisson case \cite{GR}, but in the relativistic case such compactness results 
have so far been absent in the literature. 
An essential part of the present paper is about obtaining such compactness results for minimizing sequences. 
Moreover, we provide arguments which give strong support that a class of minimizers are static solutions; 
we are however not able to fully prove this. Furthermore, our existence proof leads us to a new criterion 
which we conjecture to be sufficient for stability. 

We now specify in detail the problem we solve and we also include some more comments on the Newtonian case, 
since there are some fundamental differences between the relativistic and the non-relativistic systems. 
Let the Casimir function be defined as 
\[ \chi(f)=\frac{k}{k+1}\,f^{1+\frac{1}{k}} \] 
for $k\in ]0, 2]$. The functional we are dealing with is
\begin{equation}\label{calH1} 
   {\cal D}(f)=\int_{\R^3}\int_{\R^3} e^{\lambda}\,(\chi(f)-\alpha f)\,dx\,dv,
\end{equation} 
where $\alpha>0$, $\lambda=\lambda_f$ is given by
\[ e^{-2\lambda}=1-\frac{2m}{r} \] 
and $m=m_f$ is the mass function defined as
\[ m(r)=\int_{|x|\le r}\int_{\R^3}\sqrt{1+|v|^2}\,f\,dx\,dv. \]
We refer to this functional as the particle number-Casimir functional (introduced by Wolansky \cite{Wol} 
and also used in \cite{AK1}), since $\int_{\R^3}\int_{\R^3} e^{\lambda}\,f\,dx\,dv$ 
represents the particle number. We point out that $\mathcal{D}(f)$ is a conserved quantity along solutions $f=f(t)$ 
of the full Einstein-Vlasov system. From now on we fix $\alpha=1$, 
but let us nevertheless explain the reason for including a general $\alpha$ in (\ref{calH1}).  
The goal is to show that a minimizer of $\mathcal{D}$ over a certain set ${\cal A}$ of admissible functions 
is a static solution that has the form
\begin{equation}\label{fastform1} 
   f_0(r, w, \ell)
   =\Big(1-Ce^{\mu_0(r)}\sqrt{1+w^2+\ell^2/r^2}\Big)_+^k, 
\end{equation}
where $r=|x|$, $w=\frac{x\cdot v}{r}$ and $\ell=|x\wedge v|$ are the spherical coordinates 
and $e^{\mu_0(r)}$ is one of the metric coefficients. 
A family of static solutions of this form can be parametrized by the central red shift \cite{AR,GRS}. 
If the mass and the radius of the spatial support of the solutions for such a sequence is plotted (the `mass-radius diagram'), 
then a spiral is obtained \cite{AR}. A solution of the form (\ref{fastform1}) 
can thus not have an arbitrary large mass. However, if $f_0$ 
is a static solution with induced mass $M$, then another static solution can be generated 
by applying the scaling $\tilde{f}(x, v)=\gamma^2 f_0(\gamma x, v)$, see Lemma \ref{scalem} ahead. 
The induced mass of $\tilde{f}$ is $M/\gamma$. Hence, if we want to find static solutions with large masses as minimizers, 
then the corresponding particle number-Casimir functional (\ref{calH1}) needs to be modified, i.e., 
one needs to have a different value than $\alpha=1$. This will result in a slightly different form of $f_0$ which depends 
on $\alpha$, cf.~\cite{Wol}. For minimizing the functional ${\cal D}$, 
it is however sufficient to deal with just one $\alpha$, which we take to be $\alpha=1$. 

We will consider functions $f=f(x, v)\ge 0$ of $x, v\in\R^3$ that are spherically symmetric, 
and our goal is to minimize ${\cal D}$ 
over an appropriate class of functions ${\cal A}$. For $f\in {\cal A}$ 
we will require the mass constraint 
\begin{equation}\label{emdef} 
   \int_{\R^3}\int_{\R^3}\sqrt{1+|v|^2}\,f\,dx\,dv=M.
\end{equation} 
We note that also $M$ is a conserved quantity of the Einstein-Vlasov system.  
This constraint will not be sufficient to enforce the mass function $m=m_f$ to have 
$2m/r$ stay away from $1$. Without a restriction on $m/r$, the functional can be shown to be unbounded from below. 
This is in stark contrast to the Newtonian case, where the corresponding functional is bounded from below 
and a \textit{global} minimizer exists. 
For this reason we will need to add another condition, 
which we choose to be of the form $\rho(r)=\rho_f(r)\le\sigma_0$ for $r\in [0, \infty[$, where 
\begin{equation}\label{emden} 
   \rho(x)=\int_{\R^3}\sqrt{1+|v|^2}\,f(x, v)\,dv
\end{equation} 
denotes the density induced by $f$, and $\rho(x)=\rho(r)$ by abuse of notation. 
Here another important technical difference 
to the Vlasov-Poisson case arises: For Einstein-Vlasov (\ref{emden}) and (\ref{emdef}) 
only gives $\rho\in L^1_x(\R^3)$, whereas the density $\rho_{{\rm VP}}(x)=\int_{\R^3} f(x, v)\,dv$ 
for Vlasov-Poisson is know to satisfy $\rho_{{\rm VP}}\in L^{5/3}_x(\R^3)$.  

To summarize, we consider the convex set ${\cal A}$ given by 
\begin{eqnarray}\label{calAdef} 
   {\cal A}={\cal A}_{M,\sigma_0} & = & \Big\{f=f(x, v)\ge 0\,\,\mbox{is measurable and spherically symmetric},
   \nonumber
   \\ & & \hspace{1em} \int_{\R^3}\int_{\R^3}\sqrt{1+|v|^2}\,f\,dx\,dv=M,
   \,\,\rho(r)\le\sigma_0\,\,\mbox{for}\,\,r\in [0, \infty[\Big\}, 
\end{eqnarray}
where $M>0$ and $\beta\in ]0, \frac{1}{2}[$ are fixed, and for technical reasons 
\begin{equation}\label{siggi}  
   0<\sigma_0\le\min\Big\{1, \frac{3\beta^3}{4\pi M^2}\Big\}. 
\end{equation}
It is shown in Lemma \ref{2mrlem} below that for $f\in\mathcal{A}$ it holds that $m/r\leq\beta$. 
Let us point out that the formulation of the variational problem by Wolansky in \cite{Wol} 
is different from ours. In fact, Wolansky eventually considers another functional 
(defined in terms of a suitable Lagrangian). However, in both formulations a restriction on $m/r$ is necessary. 
This restriction introduces additional `boundaries' of the set ${\cal A}$, and the minimization problem 
is turned into an obstacle problem, which is one reason for the flaws in \cite{Wol}. 

We want to minimize $\mathcal{D}$ over the set $\mathcal{A}$ and we thus define 
\begin{equation}\label{Idef} 
   I=\inf_{f\in {\cal A}} {\cal D}(f).
\end{equation} 
It is not difficult to establish that $I<0$; see Lemma \ref{Ikle} ahead. 
The main purpose of this paper is to show that $I$ is attained, 
in the sense that there exists a function $f_0\in {\cal A}$ such that ${\cal D}(f_0)=I$. 
This will be achieved under a further assumption, 
which we call the Casimir-binding-energy condition (CBEC). 

\begin{definition}\label{sbec}  
We say that the Casimir-binding-energy condition (CBEC) holds, if 
\[ I=\inf_{f\in {\cal A}} {\cal D}(f)<-M. \]
\end{definition}
Analytically we are able to show that (CBEC) holds for all $k\in ]0, 2]$, cf.~Lemma \ref{Ikle}. 
The reason for the terminology `Casimir-binding-energy' condition 
is due to the definition of the binding energy as
\[ E_b(f)=M_0-M=\int_{\R^3}\int_{\R^3} e^{\lambda}\,f\,dx\,dv-\int_{\R^3}\int_{\R^3} \sqrt{1+|v|^2}f\,dx\,dv. \]
As mentioned above, $M_0$ is a conserved quantity called the particle number, or rest mass. 
We thus find it natural to call the quantity
\[ E_{Cb}(f)=-\mathcal{D}(f)-M=\int_{\R^3}\int_{\R^3} e^{\lambda}\,(f-\chi(f))\,dx\,dv
   -\int_{\R^3}\int_{\R^3} \sqrt{1+|v|^2}f\,dx\,dv \]
the Casimir-binding energy. 
We note that both the binding energy and the Casimir-binding energy 
are conserved quantities. Moreover, due to $f\ge 0$ we have $E_b(f)\geq E_{Cb}(f)$, so that the (CBEC) condition 
implies that a minimizer $f_0$ (which we show that exists) has the property that $E_b(f_0)>0$. 
It is known that perturbations of static solutions for which $E_b(f)>0$ have better stability properties 
than in the case where $E_b(f)<0$, although the condition $E_b(f)>0$ is certainly not sufficient for stability \cite{AR,GRS}. 
In the literature, the strong binding energy hypothesis \cite{ZP} claims that the first local maximum 
of the binding energy along a sequence of static solutions signals the onset of instability. 
This claim has recently been investigated carefully using numerical simulations, and counterexamples were found, 
see \cite{GRS}. There is however strong support that the \textit{weak} binding energy hypothesis holds, which says 
that static solutions up to the first local maximum of the binding energy along a sequence are stable. 
To decide whether or not a steady state has the property that it is situated before the first local maximum 
of the binding energy along a sequence requires that such a sequence of static solutions can be generated. 
We propose a more direct criterion that can be straightforwardly checked, 
without the need of having a sequence of static solutions. The class of perturbations that we address have the property that they leave the total mass unchanged, and the perturbed density function $f$ off a static solution $f_0$ is also required to be non-negative. We call this class mass-preserving perturbations. 
\bigskip

\textbf{Conjecture: }Let $f_0$ be an isotropic static solution of the form (\ref{fastform1}). 
If $\mathcal{D}(f_0)<-M$, where $M$ is the induced mass of $f_0$, then $f_0$ is stable for mass-preserving perturbations. 
\bigskip

In Section \ref{CBECdisc_sect} a more detailed discussion of the (CBEC) is given and some comments about the stated conjecture are also provided. 

Let us end this introduction by commenting about the relation between stability and the existence of a minimizer. 
As mentioned before, for the Vlasov-Poisson system stability is a consequence 
of the fact that a static solution is obtained as a minimizer of an energy-Casimir functional over a certain set 
(which we simply denote by $\mathcal{A}_{{\rm VP}}$ without being more specific). 
A perturbation $f$ of a static solution $f_0$ will stay in the set $\mathcal{A}_{{\rm VP}}$ as it evolves, 
and this fact is crucial for deducing stability. In the relativistic case, a perturbation of a static solution 
will not necessarily remain in the set $\mathcal{A}$ due to the constraint $\rho\leq\sigma_0$. 
Hence, if a static solution is a minimizer of the particle number-Casimir functional, 
it is necessary to control the evolution of a perturbation of the static solution 
in a much more precise fashion to make a conclusion about stability. 
This is highly non-trivial and again points at the fundamental difference between the Newtonian and the relativistic case. 
Nevertheless, this provides us with another method to approach the non-linear stability problem. 

This paper is organized as follows. In the next section the main results and some ideas of the proofs 
are outlined. A new method to obtain minimizing sequences that decay in $v$-space 
is presented in Section \ref{rearr_sect}, and it relies on a rearrangement argument. 
This property is crucial for deriving further features of the limiting function $f_0$ of the minimizing sequence, 
see Section \ref{limprop_sect}. By using the Casimir-binding-energy condition (CBEC), 
the compact support of $f_0$ in $x$-space can be proved in Section \ref{sec-spatial-part}. In addition, 
the existence of a minimizer $f_0\in \mathcal{A}$ is shown in this section. A variational inequality satisfied by $f_0$ 
is derived in Section \ref{furthprop_sect}. As a consequence of this inequality we obtain that the minimizer 
has compact support also in $v$-space. In Section \ref{isminstatic_sect}, a connection to the results in \cite{AK1} 
is made by proving that the minimizer necessarily takes a certain form, of which a static solution is a special case. 
A particularly transparent way to state the variational inequality is also derived, which clarifies under which conditions a minimizer 
is a static solution. An argument is also given which gives strong support to the belief 
that the minimizer is a static solution, provided that the mass $M$ is sufficiently small. 
A discussion of the Casimir-binding-energy condition (CBEC) is contained in Section \ref{CBECdisc_sect}, as well as a discussion of the conjecture on stability that we propose. 
Section \ref{tech_sect} has some technical results,  
in particular, it is shown that the `Casimir part' of the functional $\mathcal{D}$ is convex 
for $k\in ]0, 2]$.

%%%%%%%%%%%%%%%%%%%%%%%%%%%%%%%%%%%%%%%%%%%%%%%%%%%%%%%%%%%%%%%%%%%%%%%%%%%%%%%    

\section{Main results and ideas of proofs}
\setcounter{equation}{0}

The minimization problem (\ref{Idef}), (\ref{calAdef}), (\ref{calH1}) 
is quite non-standard. Certainly a minimizing sequence $(f_n)\subset {\cal A}$ 
will be bounded in $L^{1+1/k}(\R^3\times\R^3)$, and we can extract a subsequence 
(that is not relabeled) such that
\begin{equation}\label{kingst} 
   f_n\rightharpoonup f_0\,\,\mbox{in}\,\,L^{1+1/k}(\R^3\times\R^3)
\end{equation} 
as $n\to\infty$ for an appropriate limit function $f_0\in L^{1+1/k}(\R^3\times\R^3)$. 
However, it is far from obvious why we should have $f_0\in {\cal A}$ 
and ${\cal D}(f_0)=I$. There is also no direct concentration compactness principle 
available that could be of help at this point. 

In this paper we will first introduce a new technique that improves minimizing sequence 
by mass rearrangement in $v$-space. This procedure should be thought of 
as two `machines', with the feature that you stick in a function with certain properties, 
then turn the crank, and out comes another function that has better properties. 
The first machine modifies $f$ into $\hat{f}$ and adds the condition $\hat{f}(x, v)\le 1$, 
and in addition we have $\rho_{\hat{f}}=\rho_f$ as well as ${\cal D}(\hat{f})\le {\cal D}(f)$. 
Then the second machine takes $\hat{f}$ into some $\tilde{f}$ 
so that $\rho_{\tilde{f}}=\rho_{\hat{f}}$ and ${\cal D}(\tilde{f})\le {\cal D}(\hat{f})$ 
as well as $\tilde{f}(x, v)\le 1$. The newly gained property of $\tilde{f}$ is 
\begin{equation}\label{sterling} 
   \int_{\R^3}\int_{|v|\ge P+1}\sqrt{1+|v|^2}\,\tilde{f}\,dx\,dv\le\frac{2}{P^{1/4}},
   \quad P\ge P_0+1,
\end{equation} 
for a certain $P_0>0$ depending only on $k$, $M$ and $\beta$. 
If we apply this procedure along the minimizing sequence, 
we thus may suppose that additionally $f_n(x, v)\le 1$ holds, as well as (\ref{sterling}) for all $f_n$ 
(with uniform constants). It turns out that (\ref{sterling}) is sufficient, 
together with the weak convergence (\ref{kingst}), to prove the following facts 
about the limit function $f_0\in L^{1+1/k}(\R^3\times\R^3)$.

\begin{lemma}\label{gaum} The following assertions hold: 
\begin{itemize}
\item $0\le f_0\le 1$; 
\item $f_0$ is spherically symmetric; 
\item $\int_{\R^3}\int_{|v|\ge P+1}\sqrt{1+|v|^2}\,f_0\,dx\,dv
\le\frac{2}{P^{1/4}},\quad P\ge P_0+1$; 
\item $m_{f_n}(r)\to m_{f_0}(r)$ as $n\to\infty$ for $r\in [0, \infty[$, 
in particular $m_{f_0}(r)/r\le\beta$ for $r\in ]0, \infty[$; 
\item $\rho_{f_0}(r)\le\sigma_0$ for $r\in [0, \infty[$; 
\item $\lambda_{f_n}(r)\to\lambda_{f_0}(r)$ as $n\to\infty$ for $r\in [0, \infty[$. 
\end{itemize} 
\end{lemma}

While Lemma \ref{gaum} is satisfying, and (\ref{sterling}) is very useful to control the loss of mass in $v$-space, 
we need an additional input to make sure that also no mass is lost in $x$-space, 
in the sense that $\int_{\R^3}\int_{\R^3}\sqrt{1+|v|^2}\,f_0\,dx\,dv=M$ 
for the limit function $f_0$. Observe that we already know that 
$\int_{\R^3}\int_{\R^3}\sqrt{1+|v|^2}\,f_0\,dx\,dv
=\int_{\R^3}\int_{|v|\le P+1}\sqrt{1+|v|^2}\,f_0\,dx\,dv+{\cal O}(P^{-1/4})$, 
which serves as a substitute for compact support in $v$. 
To close the argument, we will be able to determine the potential loss of mass 
in $x$-space by comparison to a function $h\in {\cal A}$ such that ${\cal D}(h)<0$ 
and $|D(h)|>M$; such a function $h$ must exist due to (CBEC).  

Lastly, it is also possible to prove that in fact $f_0$ must have compact support, 
both in $x$ and in $v$. Hence we arrive at the following main result of this paper. 

\begin{theorem}\label{mainthm} 
Let $k\in ]0, 2]$. Then (CBEC) is verified 
and there is a function $f_0\in L^{1+1/k}(\R^3\times\R^3)$ 
that is a minimizer, i.e., we have $f_0\in {\cal A}$ and ${\cal D}(f_0)=I$. 
Furthermore, $f_0\le 1$, and there are $R_0, S_0>0$ such that $f_0(x, v)=0$ for a.e.~$|x|\ge R_0$, $v\in\R^3$ 
and for a.e.~$x\in\R^3$, $|v|\ge S_0$. 
\end{theorem}

Up to the compact support property in $v$-space, 
the proof of this theorem is completed in Section \ref{sec-spatial-part}, cf.~Theorem \ref{f0exists}. 
To establish that the minimizer has compact support in $v$-space we make use of the variational inequality 
satisfied by $f_0$, as outlined in Section \ref{furthprop_sect}. The compact support property in $v$-space 
is proven in Corollary \ref{vsuppcpt}. 
\smallskip

The important question whether or not the minimizer from Theorem \ref{mainthm} 
is a static solution of the Einstein-Vlasov system is considered in Section \ref{isminstatic_sect}. 
By making use of the variational inequality satisfied by $f_0$ and the results in \cite{AK1}, we will find that 
\[ f_0(r, w, \ell)
   =\Big(1-\tilde{\eps}(r)\sqrt{1+w^2+\ell^2/r^2}\Big)_+^k \] 
for the function $\tilde{\eps}(r)=(G')^{-1}(-4\pi\kappa\rho_0(r))$, 
where $\kappa=\frac{k+1}{16\pi^2}$ and 
\[
G(\eps)=\int_0^\infty\xi^2\,{(1-\eps\sqrt{1+\xi^2})}_+^{k+1}\,d\xi, \textrm{for } \eps\in\R.
\] 
This information can be used to rewrite the variational inequality as
\begin{equation}\label{varslim}
   \int_0^\infty\frac{d}{dr}\Big(e^{\lambda_0+\mu_0}(m_g-m_0)\Big)\,U\,dr\ge 0.
\end{equation}
Here $g\in\mathcal{A}$ is free to be taken such that 
${\rm supp}\,g\subset {\rm supp}\,f_0$, and $\mu_0$ and $\lambda_0$ 
are the induced metric functions for the minimizer $f_0$. In addition, 
\[ U(r)=-\left\{\begin{array}{c@{\quad:\quad}l} e^{-\mu_0(r)}\tilde{\eps}(r) & r\in [0, R_0] 
   \\ 0 & r\in ]R_0, \infty[\end{array}\right. , \] 
where $R_0$ denotes the end of the $r$-support of $f_0$.   
It turns out that if $U(r)={\rm const.}$, then $f_0$ is a static solution. Hence, if the inequality in (\ref{varslim}) 
could be turned into an equality, the claim would follow from integration by parts, if $U$ was differentiable. 
If the minimizer stays away from the boundary specified in the set $\mathcal{A}$, i.e., 
if the induced density satisfies $\rho_0(r)<\sigma_0$, 
then all of this can be made rigorous. We give an argument that for small mass $M$, the induced density $\rho_0$ 
must at least partially stay away from of the boundary. 
Up to such technical issues, we thus provide strong support 
that a minimizer corresponding to a small mass is a static solution of the Einstein-Vlasov system. 
 
%%%%%%%%%%%%%%%%%%%%%%%%%%%%%%%%%%%%%%%%%%%%%%%%%%%%%%%%%%%%%%%%%%%%%%%%%%%%%%%    

\section{Rearrangement in $v$-space}
\label{rearr_sect} 
\setcounter{equation}{0}

Throughout this section we fix $M>0$ and $\beta\in ]0, \frac{1}{2}[$. Let 
\begin{eqnarray}\label{tilda} 
   \tilde{{\cal A}}_{M,\beta} & = & \Big\{f=f(x, v)\ge 0\,\,\mbox{is measurable and spherically symmetric},
   \nonumber
   \\ & & \hspace{1em} \int_{\R^3}\int_{\R^3}\sqrt{1+|v|^2}\,f\,dx\,dv=M,
   \,\,\frac{m(r)}{r}\le\beta\,\,\mbox{for}\,\,r\in ]0, \infty[,
   \,\,\rho(r)\le 1\,\,\mbox{for}\,\,r\in [0, \infty[\Big\}.
   \nonumber \\ & & 
\end{eqnarray}
Since it will be needed below, we remark that the function 
\begin{equation}\label{sith} 
   \phi(s)=\chi(s)-s=\frac{k}{k+1}\,s^{1+1/k}-s,\quad s\in [0, \infty[,
\end{equation}  
has its zeros at $s=0$ and $s=(\frac{k+1}{k})^k$, attains its negative 
minimal value $\phi(1)=-\frac{1}{k+1}$ at $s=1$ 
and decreases in $[0, 1]$. 
\medskip 

We improve a given $f\in\tilde{{\cal A}}_{M,\beta}$ 
(to be thought of as a member of a minimizing sequence 
for ${\cal D}$) in several steps. 

\begin{lemma}\label{jj+1} 
Let $f=f(x, v)\ge 0$ be spherically symmetric 
and so that $f\in\tilde{{\cal A}}_{M,\beta}$, i.e., 
\begin{itemize} 
\item[(a)] $\int_{\R^3}\int_{\R^3}\sqrt{1+|v|^2}\,f\,dx\,dv=M$;  
\item[(b)] $m(r)/r\le\beta$ for $r\in ]0, \infty[$; 
\item[(c)] $\rho(r)\le 1$ for $r\in [0, \infty[$. 
\end{itemize} 
Then there is a spherically symmetric function 
$\hat{f}=\hat{f}(x, v)\ge 0$ such that 
\begin{itemize}
\item $\rho_{\hat{f}}=\rho_f$, in particular (a)-(c) holds for $\hat{f}$ 
and $\hat{f}\in\tilde{{\cal A}}_{M,\beta}$; 
\item ${\cal D}(\hat{f})\le {\cal D}(f)$;
\item $\hat{f}(x, v)\le 1$ for $(x, v)\in\R^3\times\R^3$.    
\end{itemize} 
\end{lemma} 
{\bf Proof\,:} To construct $\hat{f}$, we define the sets 
\[ E(x)=\{v\in\R^3: f(x, v)>1\} \] 
for $x\in\R^3$, and also 
\[ \rho_D(x)=\int_{E(x)}\sqrt{1+|v|^2}\,(f(x, v)-1)\,dv. \]
Then $0\le\rho_D(x)\le\rho(x)\le 1$. Let $H=\{v\in\R^3: a<|v|<b\}$ 
for $a=\sqrt[3]{\frac{1}{4\pi}}$ and $b=\sqrt[3]{\frac{25}{4\pi}}$, 
and moreover 
\[ H_<(x)=\Big\{v\in H: f(x, v)<\frac{1}{2}\Big\},
   \quad H_>(x)=\Big\{v\in H: f(x, v)\ge\frac{1}{2}\Big\} \] 
for $x\in\R^3$. We have $|H|=(4\pi/3)(b^3-a^3)=8$ for the Lebesgue measure of $H$, 
and since $H=H_<(x)\cup H_>(x)$ is a disjoint union, also $|H_<(x)|+|H_>(x)|=8$.  
We claim that $|H_>(x)|<4$. In the opposite case $|H_>(x)|\ge 4$ 
we would obtain 
\[ 1\ge\rho(x)\ge\int_{H_>(x)}\sqrt{1+|v|^2}\,f(x, v)\,dv
   \ge\frac{1}{2}\int_{H_>(x)}\sqrt{1+|v|^2}\,dv
   \ge\frac{1}{2}\,|H_>(x)|\ge 2, \] 
which is a contradiction. From $|H_>(x)|<4$ it follows that 
$|H_<(x)|>4$. Noting that $E(x)$ and $H_<(x)$ 
are disjoint, we define  
\begin{equation}\label{hatfdef} 
   \hat{f}(x, v)=\left\{\begin{array}{c@{\quad:\quad}c}
   1 & v\in E(x) 
   \\[1ex]
   \frac{\rho_D(x)}{\sqrt{1+|v|^2}\,|H_<(x)|}+f(x, v) & v\in H_<(x) 
   \\[1ex] 
   f(x, v) & v\in\R^3\setminus (E(x)\cup H_<(x)) 
   \end{array}\right. .
\end{equation} 
To see that $\rho_{\hat{f}}=\rho_f$, we calculate 
\begin{eqnarray*} 
   \rho_f(x) & = & \int_{\R^3}\sqrt{1+|v|^2}\,f(x, v)\,dv
   \\ & = & \int_{E(x)}\sqrt{1+|v|^2}\,f(x, v)\,dv
   +\int_{H_<(x)}\sqrt{1+|v|^2}\,f(x, v)\,dv	
   \\ & & +\,\int_{\R^3\setminus (E(x)\cup H_<(x))}\sqrt{1+|v|^2}\,f(x, v)\,dv
   \\ & = & \int_{E(x)}\sqrt{1+|v|^2}\,(f(x, v)-1)\,dv
   +\int_{E(x)}\sqrt{1+|v|^2}\,dv
   \\ & & +\,\int_{H_<(x)}\sqrt{1+|v|^2}\,f(x, v)\,dv	
   +\int_{\R^3\setminus (E(x)\cup H_<(x))}\sqrt{1+|v|^2}\,f(x, v)\,dv
   \\ & = & \rho_D(x)
   +\int_{E(x)}\sqrt{1+|v|^2}\,dv
   \\ & & +\,\int_{H_<(x)}\sqrt{1+|v|^2}\,f(x, v)\,dv	
   +\int_{\R^3\setminus (E(x)\cup H_<(x))}\sqrt{1+|v|^2}\,f(x, v)\,dv. 	
\end{eqnarray*} 
On the other hand, 
\begin{eqnarray*} 
   \rho_{\hat{f}}(x) & = & \int_{\R^3}\sqrt{1+|v|^2}\,\hat{f}(x, v)\,dv
   \\ & = & \int_{E(x)}\sqrt{1+|v|^2}\,dv
   +\int_{H_<(x)}\sqrt{1+|v|^2}\,\Big(\frac{\rho_D(x)}{\sqrt{1+|v|^2}\,|H_<(x)|}+f(x, v)\Big)\,dv
   \\ & & +\,\int_{\R^3\setminus (E(x)\cup H_<(x))}\sqrt{1+|v|^2}\,f(x, v)\,dv
   \\ & = & \int_{E(x)}\sqrt{1+|v|^2}\,dv+\rho_D(x)
   +\int_{H_<(x)}\sqrt{1+|v|^2}\,f(x, v)\,dv
   \\ & & +\,\int_{\R^3\setminus (E(x)\cup H_<(x))}\sqrt{1+|v|^2}\,f(x, v)\,dv,  	
\end{eqnarray*} 
so that both expressions agree. Next we argue that 
${\cal D}(\hat{f})\le {\cal D}(f)$ holds. Since $m_{\hat{f}}=m_f$, 
and hence $\lambda_{\hat{f}}=\lambda_f$, we need to establish the relation 
\[ \int_{\R^3}\int_{\R^3} e^{\lambda_f}\,(\chi(\hat{f})-\hat{f})\,dx\,dv
   \le\int_{\R^3}\int_{\R^3} e^{\lambda_f}\,(\chi(f)-f)\,dx\,dv. \] 
Due to (\ref{hatfdef}), this comes down to showing that 
\begin{equation}\label{rba} 
   \int_{\R^3} dx\,e^{\lambda_f} \int_{E(x)\cup H_<(x)} dv\,(\chi(\hat{f})-\hat{f})
   \le\int_{\R^3} dx\,e^{\lambda_f} \int_{E(x)\cup H_<(x)} dv\,(\chi(f)-f).
\end{equation}     
By definition, if $v\in H_<(x)$, then $f(x, v)<1/2$ and also 
$\delta(x, v)=\frac{\rho_D(x)}{\sqrt{1+|v|^2}\,|H_<(x)|}\le\frac{1}{|H_<(x)|}<1/4$. 
Since $1/4+1/2\le 1$, we deduce that $\phi(\hat{f}(x, v))=\phi(\delta(x, v)+f(x, v))
\le\phi(f(x, v))$. Furthermore, if $v\in E(x)$, then $\phi(\hat{f}(x, v))=\phi(1)
\le\phi(f(x, v))$, which implies (\ref{rba}). Due to $1/4+1/2\le 1$ 
and $f(x, v)\le 1$ for $v\not\in E(x)$, we also have $\hat{f}(x, v)\le 1$ 
for $(x, v)\in\R^3\times\R^3$. 
{\hfill$\Box$}\bigskip

Define $P_0\ge 10$ by 
\begin{equation}\label{P0cond}
   P_0=\max\Big\{10, \Big(\frac{1+M}{4\pi}\Big)^{4/3}, 
   \frac{64(k+1)^2}{1-2\beta}, 
   \frac{256 (k+1)^2 M^4}{(1-2\beta)^2}\Big\}. 
\end{equation} 

\begin{lemma}\label{compactsupp} 
Let $f=f(x, v)\ge 0$ be spherically symmetric 
and so that $f\in\tilde{{\cal A}}_{M,\beta}$, i.e., 
\begin{itemize} 
\item[(a)] $\int_{\R^3}\int_{\R^3}\sqrt{1+|v|^2}\,f\,dx\,dv=M$;  
\item[(b)] $m(r)/r\le\beta$ for $r\in ]0, \infty[$; 
\item[(c)] $\rho(r)\le 1$ for $r\in [0, \infty[$. 
\end{itemize} 
In addition, suppose that 
\begin{itemize} 
\item[(d)] $f(x, v)\le 1$ for $(x, v)\in\R^3\times\R^3$; 
\item[(e)] there is a $P\ge P_0$ such that 
\[ P^{1/4}\int_{\R^3}\int_{|v|\ge P}\sqrt{1+|v|^2}\,f\,dx\,dv\ge 1. \] 
\end{itemize} 
Then there is a spherically symmetric function 
$\tilde{f}=\tilde{f}(x, v)\ge 0$ such that 
\begin{itemize}
\item $\rho_{\tilde{f}}=\rho_f$, in particular (a)-(c) holds for $\tilde{f}$ 
and $\tilde{f}\in\tilde{{\cal A}}_{M,\beta}$; 
\item ${\cal D}(\tilde{f})\le {\cal D}(f)$;
\item ${\rm supp}(\tilde{f})\subset\{(x, v)\in\R^3\times\R^3: |v|\le P+1\}$;  
\item $\tilde{f}(x, v)\le 1$ for $(x, v)\in\R^3\times\R^3$;
\item we have 
\begin{equation}\label{algu} 
   P^{1/4}\int_{\R^3}\int_{P\le |v|\le P+1}\sqrt{1+|v|^2}\,\tilde{f}\,dx\,dv
   \le 1. 
\end{equation} 
\end{itemize} 
\end{lemma} 
{\bf Proof\,:} We decompose 
\begin{eqnarray*} 
   M & = & \int_{\R^3}\int_{\R^3}\sqrt{1+|v|^2}\,f\,dx\,dv
   \\ & = & \int_{\R^3}\int_{|v|\le P}\sqrt{1+|v|^2}\,f\,dx\,dv
   +\int_{\R^3}\int_{|v|\ge P}\sqrt{1+|v|^2}\,f\,dx\,dv
   \\ & =: & M_< + M_>. 
\end{eqnarray*} 
Note that we have $M_> >0$, since otherwise (e) could not hold. 
Next we introduce the numbers 
\[ N=\bigg[\frac{(1+M)P}{M_>}\bigg]\in\N,\quad\delta=\frac{1}{N},
   \quad q=\sqrt{P}, \] 
along with the strips 
\[ S_j=\Big\{(x, v)\in\R^3\times\R^3: \frac{q}{2}+j\delta q
   \le |v|\le\frac{q}{2}+(j+1)\delta q\Big\}, 
   \quad j=0, \ldots, N-1; \]
then 
\begin{equation}\label{unions} 
   \bigcup_{j=0}^{N-1} S_j=\Big\{(x, v)\in\R^3\times\R^3: 
   \frac{q}{2}\le |v|\le\frac{3q}{2}\Big\}.
\end{equation}  
Observe that $M_< \le M$ yields $\frac{(1+M)P}{M_>}\ge P\ge P_0\ge 10$, 
so at least $N\ge 9$. Furthermore, 
\begin{equation}\label{Nlow} 
   N=\bigg[\frac{(1+M)P}{M_>}\bigg]\ge\frac{(1+M)P}{M_>}-1
   \ge\frac{(1+M)P}{2M_>}
\end{equation}  
is found. We also split up 
\begin{eqnarray*} 
   {\cal D}(f) 
   & = & \int_{\R^3}\int_{\R^3} e^{\lambda}\,(\chi(f)-f)\,dx\,dv
   \\ & = & \int_{\R^3}\int_{|v|\le P} e^{\lambda}\,(\chi(f)-f)\,dx\,dv
   +\int_{\R^3}\int_{|v|\ge P} e^{\lambda}\,(\chi(f)-f)\,dx\,dv
   \\ & =: & {\cal D}_< + {\cal D}_> . 
\end{eqnarray*}
Since $f(x, v)\le 1$ pointwise, 
\begin{equation}\label{giwe}
   \chi(f(x, v))-f(x, v)=\phi(f(x, v))\le 0
\end{equation} 
pointwise; recall (\ref{sith}). In particular, ${\cal D}_< \le 0$ and ${\cal D}_> \le 0$. 
Due to $m(r)/r\le\beta$ we obtain $e^{-2\lambda}=1-2m/r\ge 1-2\beta$, so that
\begin{equation}\label{cbetabd}
   1\le e^\lambda\le\frac{1}{\sqrt{1-2\beta}}=c_\beta.
\end{equation}
It follows that 
\begin{eqnarray}\label{rbg2}  
   {\cal D}_> & = & \int_{\R^3}\int_{|v|\ge P} e^{\lambda}\,(\chi(f)-f)\,dx\,dv
   \ge -\int_{\R^3}\int_{|v|\ge P} e^{\lambda}\,f\,dx\,dv
   \ge -c_\beta\int_{\R^3}\int_{|v|\ge P} f\,dx\,dv
   \nonumber
   \\ & \ge & -\frac{c_\beta}{\sqrt{1+P^2}}\int_{\R^3}\int_{|v|\ge P}
   \sqrt{1+|v|^2}\,f\,dx\,dv
   =-\frac{c_\beta}{\sqrt{1+P^2}}\,M_>\,.
\end{eqnarray}
Furthermore, 
\begin{eqnarray}\label{rbg1}  
   |{\cal D}_<| & = & \int_{\R^3}\int_{|v|\le P} e^{\lambda}\,(f-\chi(f))\,dx\,dv
   \le\int_{\R^3}\int_{|v|\le P} e^{\lambda}\,f\,dx\,dv
   \le c_\beta\int_{\R^3}\int_{|v|\le P} f\,dx\,dv
   \nonumber
   \\ & \le & c_\beta M_<\le c_\beta M.  
\end{eqnarray} 
Now we assert that there is an $i\in\{0, \ldots, N-1\}$ such that 
\begin{equation}\label{cka} 
   \int\int_{S_i} e^{\lambda}\,(\chi(f)-f)\,dx\,dv\ge\frac{{\cal D}_<}{N}.
\end{equation}   
Otherwise we had 
\[ \int\int_{S_j} e^{\lambda}\,(\chi(f)-f)\,dx\,dv<\frac{{\cal D}_<}{N} \] 
for all $j\in\{0, \ldots, N-1\}$, which upon summation and on account 
of (\ref{unions}) would yield
\[ \int\int_{\frac{q}{2}\le |v|\le\frac{3q}{2}} 
   e^{\lambda}\,(\chi(f)-f)\,dx\,dv<{\cal D}_< . \] 
From $3q/2=3\sqrt{P}/2\le P$ in conjunction with (\ref{giwe}), 
it follows that 
\[ {\cal F}_< =\int\int_{|v|\le P} e^{\lambda}\,(\chi(f)-f)\,dx\,dv
   \le\int\int_{\frac{q}{2}\le |v|\le\frac{3q}{2}} 
   e^{\lambda}\,(\chi(f)-f)\,dx\,dv<{\cal D}_<, \] 
which is a contradiction and completes the proof of (\ref{cka}). 

Let 
\[ \rho_>(x)=\int_{|v|\ge P}\sqrt{1+|v|^2}\,f(x, v)\,dv,
   \quad\rho_i(x)=\int_{\bar{S}_i}\sqrt{1+|v|^2}\,f(x, v)\,dv, \] 
where 
\[ \bar{S}_i=\Big\{v\in\R^3: \frac{q}{2}+i\delta q
   \le |v|\le\frac{q}{2}+(i+1)\delta q\Big\} \] 
denotes the $v$-section of $S_i$. Then both $\rho_>$ and $\rho_i$ are radial, 
and also $\bar{S}_i\subset\{v\in\R^3: |v|\le 3q/2\}\le\{v\in\R^3: |v|\le P\}$. 
We claim: 
\begin{itemize}
\item[(A)] There is a $\xi_1\in ]P, P+1[$ such that 
\[ \int_{P\le |v|\le \xi_1}\sqrt{1+|v|^2}\,dv=1. \]   
\item[(B)] There is a $\xi_2\in ]\frac{q}{2}+i\delta q, \frac{q}{2}+(i+1)\delta q[$ 
such that 
\[ \int_{\frac{q}{2}+i\delta q\le |v|\le\xi_2}\sqrt{1+|v|^2}\,dv=1. \]   
\end{itemize} 
To establish (A), $\varphi_1(\xi)=\int_{P\le |v|\le\xi}\sqrt{1+|v|^2}\,dv$ 
for $\xi\in [P, P+1]$ has $\varphi_1(P)=0$ and 
\begin{eqnarray*} 
   \varphi_1(P+1) & = & \int_{P\le |v|\le P+1}\sqrt{1+|v|^2}\,dv
   =4\pi\int_P^{P+1} u^2\sqrt{1+u^2}\,du\ge 4\pi\int_P^{P+1} u^3\,du
   \\ & = & \pi ((P+1)^4-P^4)\ge 4\pi P>1.  
\end{eqnarray*} 
For (B), the argument is similar if we use the function 
$\varphi_2(\xi)=\int_{\frac{q}{2}+i\delta q\le |v|\le\xi}\sqrt{1+|v|^2}\,dv$
for $\xi\in [\frac{q}{2}+i\delta q, \frac{q}{2}+(i+1)\delta q]$, 
which has $\varphi_2(\frac{q}{2}+i\delta q)=0$ and 
\begin{eqnarray*} 
   \varphi_2\Big(\frac{q}{2}+(i+1)\delta q\Big) 
   & = & \int_{\frac{q}{2}+i\delta q\le |v|
   \le\frac{q}{2}+(i+1)\delta q}\sqrt{1+|v|^2}\,dv
   \\ & \ge & 4\pi\int_{\frac{q}{2}+i\delta q}^{\frac{q}{2}+(i+1)\delta q} u^3\,du
   =\pi\,q^4
   \Big((1+2(i+1)\delta)^4-(1+2i\delta)^4\Big)
   \\ & = & \pi\,q^4
   \Big(8\delta+24[(i+1)^2-i^2]\delta^2
   +32[(i+1)^3-i^3]\delta^3+16[(i+1)^4-i^4]\delta^4\Big)
   \\ & \ge & 8\pi\,q^4\delta=8\pi\,\frac{P^2}{N}
   \ge 8\pi\,\frac{P^2 M_>}{(1+M)P}
   \ge 8\pi\,\frac{P^{3/4}}{1+M}>1, 
\end{eqnarray*} 
where we have used that $P^{1/4} M_> \ge 1$ is verified by (e), 
along with the definition of $P_0$. 
Thus (A) and (B) have been established, and we can put 
\[ h_1(x, v)={\bf 1}_{[P,\,\xi_1]}(|v|)\rho_i(x),
   \quad h_2(x, v)={\bf 1}_{[\frac{q}{2}+i\delta q,\,\xi_2]}(|v|)\rho_>(x). \] 

Define $\tilde{f}$ by means of 
\[ \tilde{f}(x, v)=\left\{\begin{array}{c@{\quad:\quad}c}
   0 & |v|\ge P+1 \\[1ex]
   h_1(x, v) & P\le |v|\le P+1\\[1ex] 
   f(x, v) & (x, v)\in \{(x, v): |v|\le P\}\setminus S_i\\[1ex]  
   h_2(x, v) & (x, v)\in S_i    
   \end{array}\right. . \]
To see that $\tilde{f}$ has the desired properties, we begin by calculating 
its density $\rho_{\tilde{f}}$. Here we have, by (A) and (B) above, 
\begin{eqnarray*} 
   \rho_{\tilde{f}}(x) & = & \int_{\R^3}\sqrt{1+|v|^2}\,\tilde{f}(x, v)\,dv 
   \\ & = & \int_{P\le |v|\le P+1}\sqrt{1+|v|^2}\,h_1(x, v)\,dv 
   +\int_{\{|v|\le P\}\setminus\bar{S}_i}\sqrt{1+|v|^2}\,f(x, v)\,dv 
   \\ & & +\,\int_{\bar{S}_i}\sqrt{1+|v|^2}\,h_2(x, v)\,dv 
   \\ & = & \rho_i(x)\int_{P\le |v|\le\xi_1}\sqrt{1+|v|^2}\,dv 
   +\int_{\{|v|\le P\}\setminus\bar{S}_i}\sqrt{1+|v|^2}\,f(x, v)\,dv 
   \\ & & +\rho_>(x)\int_{\frac{q}{2}+i\delta q\le |v|\le\xi_2}
   \sqrt{1+|v|^2}\,dv 
   \\ & = & \int_{\bar{S}_i}\sqrt{1+|v|^2}\,f(x, v)\,dv 
   +\int_{\{|v|\le P\}\setminus\bar{S}_i}\sqrt{1+|v|^2}\,f(x, v)\,dv 
   \\ & & +\int_{|v|\ge P}\sqrt{1+|v|^2}\,f(x, v)\,dv
   \\ & = & \int_{\R^3}\sqrt{1+|v|^2}\,f(x, v)\,dv=\rho_f(x),  
\end{eqnarray*} 
so that in particular $m_{\tilde{f}}=m_f$ and consequently 
$\lambda_{\tilde{f}}=\lambda_f$. We can also write 
\begin{eqnarray}\label{newH}  
   {\cal D}(\tilde{f})
   & = & \int_{\R^3}\int_{\R^3} e^{\lambda_{\tilde{f}}}
   \,(\chi(\tilde{f})-\tilde{f})\,dx\,dv
   =\int_{\R^3}\int_{\R^3} e^{\lambda_f}\,(\chi(\tilde{f})-\tilde{f})\,dx\,dv
   \nonumber
   \\ & = & \int_{\R^3}\int_{P\le |v|\le P+1} 
   e^{\lambda_f}\,(\chi(h_1)-h_1)\,dx\,dv
   +\int_{\R^3}\int_{\{|v|\le P\}\setminus\bar{S}_i} e^{\lambda_f}
   \,(\chi(f)-f)\,dx\,dv
   \nonumber
   \\ & & +\int_{\R^3}\int_{\bar{S}_i} e^{\lambda_f}\,(\chi(h_2)-h_2)\,dx\,dv. 
\end{eqnarray}  
For the contribution on $\{P\le |v|\le P+1\}$ we note that $0\le\rho_i(x)\le\rho(x)\le 1$ 
yields $0\le h_1(x, v)\le 1$, and accordingly $\chi(h_1(x, v))-h_1(x, v)
=\phi(h_1(x, v))\le 0$. Thus  
\begin{equation}\label{newH1} 
   \int_{\R^3}\int_{P\le |v|\le P+1} 
   e^{\lambda_f}\,(\chi(h_1)-h_1)\,dx\,dv\le 0.
\end{equation}  
For the part on $\{|v|\le P\}\setminus\bar{S}_i$ we use (\ref{cka}) to get 
\begin{eqnarray}\label{newH2} 
   \lefteqn{\int_{\R^3}\int_{\{|v|\le P\}\setminus\bar{S}_i} e^{\lambda_f}
   \,(\chi(f)-f)\,dx\,dv}
   \nonumber
   \\ & = & \int_{\R^3}\int_{|v|\le P} e^{\lambda_f}
   \,(\chi(f)-f)\,dx\,dv
   -\int_{\R^3}\int_{\bar{S}_i} e^{\lambda_f}
   \,(\chi(f)-f)\,dx\,dv
   \nonumber
   \\ & = & {\cal D}_< -\int\int_{S_i} e^{\lambda_f}
   \,(\chi(f)-f)\,dx\,dv
   \nonumber
   \\ & \le & {\cal D}_< - \frac{{\cal D}_<}{N}. 
\end{eqnarray} 
Lastly, as for $h_1$ we obtain $0\le h_2(x, v)\le 1$, so that 
\begin{equation}\label{h2low} 
   h_2-\chi(h_2)=h_2\Big(1-\frac{k}{k+1}\,h_2^{1/k}\Big)\ge\frac{1}{k+1}\,h_2.
\end{equation}  
Hence (\ref{cbetabd}), the definition of $h_2$ and (B) above yields 
\begin{eqnarray}\label{newH3} 
   \int_{\R^3}\int_{\bar{S}_i} e^{\lambda_f}\,(h_2-\chi(h_2))\,dx\,dv
   & \ge & \frac{1}{k+1}\int_{\R^3}\int_{\bar{S}_i} e^{\lambda_f}\,h_2\,dx\,dv
   \ge\frac{1}{k+1}\int_{\R^3}\int_{\bar{S}_i} h_2\,dx\,dv
   \nonumber
   \\ & = & \frac{1}{k+1}\Big(\int_{\R^3}\rho_>(x)\,dx\Big)
   \Big(\int_{\bar{S}_i}{\bf 1}_{[\frac{q}{2}+i\delta q,\,\xi_2]}(|v|)\,dv\Big)
   \nonumber
   \\ & = & \frac{1}{k+1}\,M_>
   \Big(\int_{\frac{q}{2}+i\delta q\le |v|\le\xi_2}
   \,\frac{\sqrt{1+|v|^2}}{\sqrt{1+|v|^2}}\,dv\Big)
   \nonumber 
   \\ & \ge & \frac{1}{k+1}\,\frac{M_>}{\sqrt{1+\frac{9q^2}{4}}}
   \Big(\int_{\frac{q}{2}+i\delta q\le |v|\le\xi_2}\,\sqrt{1+|v|^2}\,dv\Big)
   \nonumber
   \\ & = & \frac{1}{k+1}\,\frac{M_>}{\sqrt{1+\frac{9q^2}{4}}}. 
\end{eqnarray} 
Now if we put (\ref{newH1}), (\ref{newH2}) and (\ref{newH3}) 
into (\ref{newH}), it follows that
\[ {\cal D}(\tilde{f})
   \le {\cal D}_< - \frac{{\cal D}_<}{N}
   -\frac{1}{k+1}\,\frac{M_>}{\sqrt{1+\frac{9q^2}{4}}}. \] 
Since $q^2=P$ and $1+9P/4\le 4P$, 
we deduce from (\ref{Nlow}), (\ref{rbg1}) and (\ref{rbg2}) that 
\begin{eqnarray*} 
   {\cal D}(\tilde{f})
   & \le & {\cal D}_< +\frac{2|{\cal D}_<| M_>}{(1+M)P}
   -\frac{1}{2}\,\frac{1}{k+1}\,\frac{M_>}{\sqrt{P}}
   \\ & = & {\cal D}_< +\frac{2|{\cal D}_<| M_>}{(1+M)P}
   -\frac{1}{4}\,\frac{1}{k+1}\,\frac{M_>}{\sqrt{P}}
   -\frac{1}{4}\,\frac{1}{k+1}\,\frac{M_>}{\sqrt{P}}
   \\ & \le & {\cal D}_< + {\cal D}_> +\frac{M_>}{\sqrt{P}}
   \bigg[\frac{2c_\beta}{\sqrt{P}}-\frac{1}{4}\,\frac{1}{k+1}\bigg]
   +\frac{M_>}{\sqrt{P}}\bigg[\frac{c_\beta\sqrt{P}}{\sqrt{1+P^2}}
   -\frac{1}{4}\,\frac{1}{k+1}\bigg]
   \\ & \le & {\cal D}_< + {\cal D}_> = {\cal D}(f), 
\end{eqnarray*} 
where we have also used the definition (\ref{P0cond}) of $P_0$. 
The asserted support property of $\tilde{f}$ is just by definition, 
and so is $\tilde{f}(x, v)\le 1$, since $h_1(x, v)\le 1$, 
$h_2(x, v)\le 1$ as well as $f(x, v)\le 1$. 
Finally to establish (\ref{algu}), we get from (A), 
and the definitions of $\rho_i$ and $\bar{S}_i$,   
\begin{eqnarray*}
   \int_{\R^3}\int_{P\le |v|\le P+1}\sqrt{1+|v|^2}\,\tilde{f}\,dx\,dv
   & = & \int_{\R^3}\int_{P\le |v|\le P+1}\sqrt{1+|v|^2}\,h_1\,dx\,dv
   \\ & = & \bigg(\int_{P\le |v|\le \xi_1}\sqrt{1+|v|^2}\,dv\bigg)
   \bigg(\int_{\R^3}\rho_i(x)\,dx\bigg) 
   \\ & = & \int_{\R^3}\int_{\bar{S}_i}\sqrt{1+|v|^2}\,f\,dx\,dv
   \\ & = & \int\int_{S_i}\sqrt{1+|v|^2}\,f\,dx\,dv. 
\end{eqnarray*} 
Recalling that $f(x, v)\le 1$, the analogous argument as is 
(\ref{h2low}) shows that $f-\chi(f)\ge\frac{1}{k+1}\,f$. 
Furthermore, on $S_i$ we have $|v|\le 3q/2$. 
By means of (\ref{cbetabd}), (\ref{cka}), (\ref{rbg1}) 
and (\ref{Nlow}) we may therefore bound 
\begin{eqnarray*}
   \lefteqn{P^{1/4}\int_{\R^3}\int_{P\le |v|\le P+1}\sqrt{1+|v|^2}\,\tilde{f}\,dx\,dv} 
   \\ & \le & (k+1) \sqrt{1+\frac{9q^2}{4}}
   \,P^{1/4}\int\int_{S_i} (f-\chi(f))\,dx\,dv 
   \\ & \le & (k+1) \sqrt{1+\frac{9q^2}{4}}
   \,P^{1/4}\int\int_{S_i} e^\lambda (f-\chi(f))\,dx\,dv 
   \\ & \le & (k+1) \sqrt{1+\frac{9q^2}{4}}
   \,P^{1/4}\,\frac{|{\cal D}_<|}{N}
   \\ & \le & (k+1) \sqrt{1+\frac{9q^2}{4}}
   \,P^{1/4}\,\frac{c_\beta M}{N}
   \\ & \le & (k+1) \sqrt{1+\frac{9P}{4}}
   \,P^{1/4}\,\frac{2c_\beta M M_>}{(1+M)P}
   \\ & \le & (k+1)\,\frac{4c_\beta M}{P^{1/4}}
   \\ & \le & 1, 
\end{eqnarray*} 
and this completes the proof. 
{\hfill$\Box$}\bigskip

\begin{cor}\label{rahner} 
Let $f=f(x, v)\ge 0$ be spherically symmetric 
and so that $f\in\tilde{{\cal A}}_{M,\beta}$, i.e., 
\begin{itemize} 
\item[(a)] $\int_{\R^3}\int_{\R^3}\sqrt{1+|v|^2}\,f\,dx\,dv=M$;  
\item[(b)] $m(r)/r\le\beta$ for $r\in ]0, \infty[$; 
\item[(c)] $\rho(r)\le 1$ for $r\in [0, \infty[$. 
\end{itemize} 
Then there is a spherically symmetric function 
$\tilde{f}=\tilde{f}(x, v)\ge 0$ such that 
\begin{itemize}
\item $\rho_{\tilde{f}}=\rho_f$, in particular (a)-(c) holds for $\tilde{f}$ 
and $\tilde{f}\in\tilde{{\cal A}}_{M,\beta}$; 
\item ${\cal D}(\tilde{f})\le {\cal D}(f)$;
\item $\tilde{f}(x, v)\le 1$ for $(x, v)\in\R^3\times\R^3$;  
\item we have 
\begin{equation}\label{allP} 
   \int_{\R^3}\int_{|v|\ge P+1}\sqrt{1+|v|^2}\,\tilde{f}\,dx\,dv\le\frac{2}{P^{1/4}},
   \quad P\ge P_0+1.
\end{equation}     
\end{itemize} 
\end{cor} 
{\bf Proof\,:} From $f$ we can pass to $\hat{f}$ 
by applying Lemma \ref{jj+1}, and hence we may suppose 
w.l.o.g.~that already $f$ has the additional property that 
\begin{equation}\label{apl1}
   f(x, v)\le 1\quad\mbox{for}\quad (x, v)\in\R^3\times\R^3. 
\end{equation} 
Now we need to distinguish a few cases. 
\underline{Case 1:} We have 
\[ P_0^{1/4}\int_{\R^3}\int_{|v|\ge P_0}\sqrt{1+|v|^2}\,f\,dx\,dv\ge 1. \]
Then (a)-(e) from Lemma \ref{compactsupp} are satisfied 
for $f$ and $P_0$. Thus there is $\tilde{f}$ such that 
$\rho_{\tilde{f}}=\rho_f$, ${\cal D}(\tilde{f})\le {\cal D}(f)$, 
${\rm supp}(\tilde{f})\subset\{(x, v)\in\R^3\times\R^3: |v|\le P_0+1\}$ 
and $\tilde{f}(x, v)\le 1$ pointwise; 
no further property of $\tilde{f}$ will be needed here.   
To prove (\ref{allP}), let $P\ge P_0+1$. If $|v|\ge P$, 
then $\tilde{f}(x, v)=0$ and hence 
\[ P^{1/4}\int_{\R^3}\int_{|v|\ge P}\sqrt{1+|v|^2}\,\tilde{f}\,dx\,dv=0. \] 
\underline{Case 2:} We have 
\[ P_0^{1/4}\int_{\R^3}\int_{|v|\ge P_0}\sqrt{1+|v|^2}\,f\,dx\,dv<1 \]
and
\[ P^{1/4}\int_{\R^3}\int_{|v|\ge P}\sqrt{1+|v|^2}\,f\,dx\,dv\le 1,
   \quad P\ge P_0, \] 
holds for $f$. Since 
\[ \int_{\R^3}\int_{|v|\ge P+1}\sqrt{1+|v|^2}\,f\,dx\,dv
   \le\int_{\R^3}\int_{|v|\ge P}\sqrt{1+|v|^2}\,f\,dx\,dv \] 
we may simply take $\tilde{f}=f$ in this case. 
\underline{Case 3:} We have 
\[ P_0^{1/4}\int_{\R^3}\int_{|v|\ge P_0}\sqrt{1+|v|^2}\,f\,dx\,dv<1 \] 
and there is a $P_1\ge P_0$ such that 
\[ P_1^{1/4}\int_{\R^3}\int_{|v|\ge P_1}\sqrt{1+|v|^2}\,f\,dx\,dv>1 \] 
is verified. Then $P_1>P_0$, and since the function 
$\varphi(P)=P^{1/4}\int_{\R^3}\int_{|v|\ge P}\sqrt{1+|v|^2}\,f\,dx\,dv$ 
is continuous on $[P_0, P_1]$, there is a first $\hat{P}\in ]P_0, P_1[$ 
such that 
\[ \hat{P}^{1/4}\int_{\R^3}\int_{|v|\ge\hat{P}}\sqrt{1+|v|^2}\,f\,dx\,dv=1. \] 
In particular, 
\begin{equation}\label{1a3} 
   P^{1/4}\int_{\R^3}\int_{|v|\ge P}\sqrt{1+|v|^2}\,f\,dx\,dv\le 1,
   \quad P\in [P_0, \hat{P}].
\end{equation}  
For $f$ and $\hat{P}$ we can apply Lemma \ref{compactsupp} 
to obtain a density $\tilde{f}$ such that 
$\rho_{\tilde{f}}=\rho_f$, ${\cal D}(\tilde{f})\le {\cal D}(f)$, 
${\rm supp}(\tilde{f})\subset\{(x, v)\in\R^3\times\R^3: |v|\le\hat{P}+1\}$ 
and $\tilde{f}(x, v)\le 1$ pointwise. In addition, we can use that 
\begin{equation}\label{alguzu} 
   \hat{P}^{1/4}\int_{\R^3}\int_{\hat{P}\le |v|\le\hat{P}+1}\sqrt{1+|v|^2}\,\tilde{f}\,dx\,dv
   \le 1. 
\end{equation}  
To show (\ref{allP}), let $P\ge P_0+1$. If $P\ge \hat{P}$ and $|v|\ge P+1$, 
then $|v|\ge\hat{P}+1$, and thus the integral in question vanishes 
by the support properties of $\tilde{f}$. Hence we may assume that 
$P_0+1\le P\le\hat{P}$ holds. 
First we will consider the case where $P+1\le\hat{P}$, 
so that we have $P_0\le P<P+1\le\hat{P}<\hat{P}+1$. We write 
\begin{eqnarray}\label{wema}  
  \lefteqn{P^{1/4}\int_{\R^3}\int_{|v|\ge P+1}\sqrt{1+|v|^2}\,\tilde{f}\,dx\,dv}
  \nonumber
  \\ & = & P^{1/4}\int_{\R^3}\int_{P+1\le |v|\le\hat{P}+1}\sqrt{1+|v|^2}\,\tilde{f}\,dx\,dv
  \nonumber
  \\ & = & P^{1/4}\int_{\R^3}\int_{P+1\le |v|\le\hat{P}}\sqrt{1+|v|^2}\,\tilde{f}\,dx\,dv
  +P^{1/4}\int_{\R^3}\int_{\hat{P}\le |v|\le\hat{P}+1}\sqrt{1+|v|^2}\,\tilde{f}\,dx\,dv
  \qquad 
\end{eqnarray} 
Now 
\[ P^{1/4}\int_{\R^3}\int_{P+1\le |v|\le\hat{P}}\sqrt{1+|v|^2}\,f\,dx\,dv
   \le P^{1/4}\int_{\R^3}\int_{|v|\ge P}\sqrt{1+|v|^2}\,f\,dx\,dv\le 1 \]  
by (\ref{1a3}), as $P\le\hat{P}$. 
The second term in (\ref{wema}) is also bounded by $1$, 
since $P^{1/4}\le\hat{P}^{1/4}$ and we may use (\ref{alguzu}). 
It remains to consider the case where $P+1\ge\hat{P}$, 
so that $P_0\le P\le\hat{P}\le P+1\le\hat{P}+1$. Then 
\begin{eqnarray*}
  P^{1/4}\int_{\R^3}\int_{|v|\ge P+1}\sqrt{1+|v|^2}\,\tilde{f}\,dx\,dv
  & = & P^{1/4}\int_{\R^3}\int_{P+1\le |v|\le\hat{P}+1}\sqrt{1+|v|^2}\,\tilde{f}\,dx\,dv
  \\ & \le & \hat{P}^{1/4}\int_{\R^3}\int_{\hat{P}\le |v|\le\hat{P}+1}\sqrt{1+|v|^2}\,\tilde{f}\,dx\,dv
  \\ & \le & 1,    
\end{eqnarray*} 
once again by (\ref{alguzu}). 
{\hfill$\Box$}\bigskip

%%%%%%%%%%%%%%%%%%%%%%%%%%%%%%%%%%%%%%%%%%%%%%%%%%%%%%%%%%%%%%%%%%%%%%%%%%%%%%%%%%%%%%%%%

\section{Properties of the limit function}
\label{limprop_sect} 
\setcounter{equation}{0}

\begin{remark}\label{AinA}{\rm If (\ref{siggi}) holds, 
then ${\cal A}_{M,\sigma_0}\subset\tilde{{\cal A}}_{M,\beta}$, 
where $\tilde{{\cal A}}_{M,\beta}$ is from (\ref{tilda}).  
To see this, $\rho(r)\le\sigma_0\le 1$ for $r\in [0, \infty[$ 
by the definition of $\sigma_0$. Recalling $\sigma_M=\frac{3}{32\pi M^2}$, 
also Lemma \ref{2mrlem} ahead shows that  
\[ \frac{2m(r)}{r}\le\bigg(\frac{\sigma_0}{\sigma_M}\bigg)^{1/3}
   \le\bigg(\frac{96\pi\beta^3M^2}{12\pi M^2}\bigg)^{1/3}=2\beta \] 
for $r\in [0, \infty[$, and hence $m(r)/r\le\beta$ for $r\in [0, \infty[$. 
}
\end{remark} 

We recall that
\[ I=\inf_{f\in {\cal A}} {\cal D}(f). \]

\begin{lemma}\label{Ikle} 
If $f\in {\cal A}$, then ${\cal D}(f)\ge -c_\beta M$ 
for $c_\beta=\frac{1}{\sqrt{1-2\beta}}$;
in particular, $I$ is well-defined. Furthermore, 
\[ I<0. \]
In addition, for every exponent $k>0$ there exists a function $h\in {\cal A}$ 
such that ${\cal D}(h)<0$ and $|{\cal D}(h)|>M$. 
\end{lemma}
{\bf Proof\,:} For $f\in {\cal A}$ we have $1\le e^\lambda\le c_\beta$, 
as in (\ref{cbetabd}), and also recall Remark \ref{AinA}. Therefore
\[ {\cal D}(f)=\int_{\R^3}\int_{\R^3} e^{\lambda}\,(\chi(f)-f)\,dx\,dv
   \ge -\int_{\R^3}\int_{\R^3} e^{\lambda}\,f\,dx\,dv\ge -c_\beta M. \]
To show that $I$ is negative, consider 
\[ \tilde{f}=A\,{\bf 1}_{[0, a]}(|x|)\,{\bf 1}_{[0, b]}(|v|) \]
for $A\in ]0, 1]$ and $a, b>0$. Then 
\begin{eqnarray*} 
   \tilde{\rho}(r) & = & \int_{\R^3}\sqrt{1+|v|^2}\,\tilde{f}\,dv
   =A\,{\bf 1}_{[0, a]}(r)\,\theta(b),
  \\ \theta(b) & = & 4\pi\int_0^b s^2\sqrt{1+s^2}\,ds
   =\frac{\pi}{2}\Big(b\sqrt{1+b^2}(1+2b^2)-\arcsin b\Big).
\end{eqnarray*} 
This gives the mass function 
\[ \tilde{m}(r)=4\pi\int_0^r s^2\tilde{\rho}(s)\,ds
   =\left\{\begin{array}{c@{\quad:\quad}c}
   \frac{4\pi}{3}\,A\,\theta(b)\,r^3 & r\le a 
   \\[1ex]
   \frac{4\pi}{3}\,A\,\theta(b)\,a^3 & r\ge a
   \end{array}\right. , \] 
and hence the condition $\tilde{f}\in {\cal A}$ is equivalent to 
\begin{equation}\label{abcond} 
   \frac{4\pi}{3}\,A\,\theta(b)\,a^3=M
   \quad\mbox{and}\quad A\,\theta(b)\le\sigma_0
   \le\min\Big\{1, \frac{3\beta^3}{4\pi M^2}\Big\}.
\end{equation}  
In addition, we obtain 
\[ e^{\tilde{\lambda}}=\frac{1}{\sqrt{1-\frac{2\tilde{m}(r)}{r}}}
   =\left\{\begin{array}{c@{\quad:\quad}c}
   \frac{1}{\sqrt{1-\frac{8\pi}{3}\,A\,\theta(b)\,r^2}} & r\le a 
   \\[2ex]
   \frac{1}{\sqrt{1-\frac{8\pi}{3}\,A\,\theta(b)\,\frac{a^3}{r}}} & r\ge a 
   \end{array}\right. . \] 
Since $\tilde{f}=A$ on its support and $\chi(A)-A=-A(1-\frac{k}{k+1}\,A^{\frac{1}{k}})<0$ 
due to $A\le 1$, this yields 
\begin{eqnarray}\label{ost}
   {\cal D}(\tilde{f}) 
   & = & \int_{\R^3}\int_{\R^3} e^{\tilde{\lambda}}
   (\chi(\tilde{f})-\tilde{f})\,dx\,dv 
   \nonumber
   \\ & = &  -A\Big(1-\frac{k}{k+1}\,A^{\frac{1}{k}}\Big)\,\frac{16\pi^2}{3}\,b^3
   \int_0^a\frac{r^2}{\sqrt{1-\frac{8\pi}{3}\,A\,\theta(b)\,r^2}}\,dr<0.  
\end{eqnarray}  
Note that here explicitly, for $r\in [0, a]$,  
\[ \frac{2\tilde{m}(r)}{r}=\frac{8\pi}{3}\,A\,\theta(b)\,r^2
   \le\frac{8\pi}{3}\,A\,\theta(b)\,a^2\le 2\beta<1 \] 
for $r\le a$, as predicted by Lemma \ref{2mrlem}, using that (\ref{abcond}) yields 
$\theta(b)=\frac{3M}{4\pi A a^3}$ and $\frac{3M}{4\pi a^3}=A\theta(b)\le\frac{3\beta^3}{4\pi M^2}$, whence 
\[ \frac{8\pi}{3}\,A\,\theta(b)\,a^2=\frac{2M}{a}\le 2\beta. \] 
The condition (\ref{abcond}) is easy to satisfy: We will fix $A=A(k)$ depending only on $k$ below. 
Then first take $b>0$ small enough such that $A\theta(b)\le\sigma_0$ holds 
and thereafter put $a=(\frac{3M}{4\pi A\theta(b)})^{1/3}$. In this case $\tilde{f}\in {\cal A}$ 
and ${\cal D}(\tilde{f})<0$ by (\ref{ost}), so that in particular $I<0$. 
\smallskip 

\noindent
For the last assertion, the integral (\ref{ost}) can be evaluated explicitly. 
It is found that, with $c=(\frac{8\pi}{3}\,A\theta(b))^{1/2}$, one has 
\[ {\cal D}(\tilde{f}) 
   = -A\Big(1-\frac{k}{k+1}\,A^{\frac{1}{k}}\Big)\,\frac{8\pi^2}{3}\,
   \frac{b^3}{c^3}\,\Big(\arcsin(ca)-ca\sqrt{1-c^2 a^2}\Big). \] 
Thus if (\ref{abcond}) holds, the condition $|{\cal D}(\tilde{f})|>M$ is equivalent to 
\begin{equation}\label{skan}  
   4\Big(1-\frac{k}{k+1}\,A^{\frac{1}{k}}\Big)\,
   \frac{\arcsin(ca)-ca\sqrt{1-c^2 a^2}}{c^3 a^3}
   >\frac{b\sqrt{1+b^2}(1+2b^2)-\arcsin b}{b^3}.  
\end{equation} 
As $x\to 0^+$, the expansions are 
\begin{eqnarray*} 
   & & \arcsin(x)-x\sqrt{1-x^2}=\frac{2}{3}\,x^3+{\cal O}(x^5), 
   \\ & & x\sqrt{1+x^2}(1+2x^2)-\arcsin x=\frac{7}{3}\,x^3+{\cal O}(x^5), 
\end{eqnarray*} 
and we already know from above that we can take $b\to 0$ 
and $a=(\frac{3M}{4\pi A\theta(b)})^{1/3}$ to make sure that $\tilde{f}\in {\cal A}$. 
From the expansions we deduce as $b\to 0$ 
\begin{eqnarray*} 
   \theta(b) & = & \frac{\pi}{2}\,\Big(\frac{7}{3}\,b^3+{\cal O}(b^5)\Big) 
   =\frac{7\pi}{6}\,b^3+{\cal O}(b^5), 
   \\ c & = & \sqrt{\frac{8\pi}{3}\,A\theta(b)} 
   =\sqrt{A}\,\frac{2\sqrt{7}\pi}{3}\,b^{3/2}\,(1+{\cal O}(b^2)), 
   \\ a & = & \Big(\frac{3M}{4\pi A\theta(b)}\Big)^{1/3}
   =A^{-1/3}\,\Big(\frac{9M}{14\pi^2}\Big)^{1/3}\,\frac{1}{b}\,(1+{\cal O}(b^2)), 
   \\ ca & = & A^{1/6}\,\frac{2\sqrt{7}\pi}{3}\,\Big(\frac{9M}{14\pi^2}\Big)^{1/3} b^{1/2}\,(1+{\cal O}(b^2)), 
\end{eqnarray*} 
so that also $ca\to 0$ as $b\to 0$. Therefore we may use the expansions to infer 
\[ \frac{\arcsin(ca)-ca\sqrt{1-c^2 a^2}}{c^3 a^3}=\frac{2}{3}+{\cal O}(b),
   \quad\frac{b\sqrt{1+b^2}(1+2b^2)-\arcsin b}{b^3}=\frac{7}{3}+{\cal O}(b^2).  \] 
Thus if 
\begin{equation}\label{17} 
   4\Big(1-\frac{k}{k+1}\,A^{\frac{1}{k}}\Big)\,\cdot\frac{2}{3}>\frac{7}{3}
\end{equation} 
holds, then (\ref{skan}) will be verified for $b>0$ small enough; 
the condition (\ref{17}) reads as 
\[ \frac{k}{k+1}\,A^{\frac{1}{k}}<\frac{1}{8}, \] 
and it can be satisfied, if $k>0$ is fixed and $A=A(k)>0$ is taken 
to be sufficiently small; for instance $A(k)=(\frac{1}{8})^k$ will do. 
{\hfill$\Box$}\bigskip

Let $(f_n)\subset {\cal A}$ be a minimizing sequence, i.e.
\begin{equation}\label{minsequ}
   \lim_{n\to\infty} {\cal D}(f_n)=I.
\end{equation}
For $M>0$ and $\beta\in ]0, \frac{1}{2}[$ given, let $P_0>0$ be determined 
by (\ref{P0cond}). Since ${\cal A}={\cal A}_{M,\sigma_0}\subset\tilde{{\cal A}}_{M,\beta}$ 
by Remark \ref{AinA}, the results from Section \ref{rearr_sect} are applicable, 
and in particular every $f_n$ satisfies (a)-(c) from Corollary \ref{rahner}. 
Hence there is $\tilde{f}_n\in\tilde{{\cal A}}_{M,\beta}$ such that $0\le\tilde{f}_n\le 1$, 
$\rho_{\tilde{f}_n}=\rho_{f_n}$, ${\cal D}(\tilde{f}_n)\le {\cal D}(f_n)$, and   
\begin{equation}\label{allPn} 
   \int_{\R^3}\int_{|v|\ge P+1}\sqrt{1+|v|^2}\,\tilde{f}_n\,dx\,dv\le\frac{2}{P^{1/4}},
   \quad P\ge P_0+1.
\end{equation}     
Since $\rho_{\tilde{f}_n}=\rho_{f_n}$, also $\tilde{f}_n\in {\cal A}$. 
Therefore ${\cal D}(\tilde{f}_n)\ge I$, and thus $\lim_{n\to\infty} {\cal H}(\tilde{f}_n)=I$, 
which shows that $(\tilde{f}_n)$ is a minimizing sequence, too. 
Thus we may simply drop the tilde and assume that 
the original sequence $(f_n)$ possesses the property (\ref{allPn}), 
and moreover $0\le f_n\le 1$ for every $n\in\N$. 

\begin{lemma}\label{1+epsbd}
$(f_n)\subset L^{1+1/k}(\R^3\times\R^3)$ is bounded.
\end{lemma}
{\bf Proof\,:} By (\ref{cbetabd}), 
\begin{eqnarray*}
   \int_{\R^3}\int_{\R^3}\chi(f_n)\,dx\,dv
   & \le & \int_{\R^3}\int_{\R^3} e^{\lambda_{f_n}}\chi(f_n)\,dx\,dv
   \\ & = & {\cal D}(f_n)+\int_{\R^3}\int_{\R^3} e^{\lambda_{f_n}}\,f_n\,dx\,dv
   \\ & \le & {\cal D}(f_n)+c_\beta M
\end{eqnarray*}
together with (\ref{minsequ}) yields the claim.
{\hfill$\Box$}\bigskip

As a consequence of Lemma \ref{1+epsbd} we may extract a subsequence (not relabeled) such that
\begin{equation}\label{weakc}
   f_n\rightharpoonup f_0\,\,\mbox{in}\,\,L^{1+1/k}(\R^3\times\R^3)
\end{equation}
as $n\to\infty$ for an appropriate limit function $f_0\in L^{1+1/k}(\R^3\times\R^3)$.
Now we are ready for the proof of Lemma \ref{gaum}, 
the statement of which we repeat for convenience. 

\begin{lemma}\label{f0props} The following assertions hold: 
\begin{itemize}
\item $0\le f_0\le 1$; 
\item $f_0$ is spherically symmetric; 
\item $\int_{\R^3}\int_{|v|\ge P+1}\sqrt{1+|v|^2}\,f_0\,dx\,dv
\le\frac{2}{P^{1/4}},\quad P\ge P_0+1$; 
\item $m_{f_n}(r)\to m_{f_0}(r)$ as $n\to\infty$ for $r\in [0, \infty[$, 
in particular $m_{f_0}(r)/r\le\beta$ for $r\in ]0, \infty[$; 
\item $\rho_{f_0}(r)\le\sigma_0$ for $r\in [0, \infty[$; 
\item $\lambda_{f_n}(r)\to\lambda_{f_0}(r)$ as $n\to\infty$ for $r\in [0, \infty[$. 
\end{itemize} 
\end{lemma}
{\bf Proof\,:} If $U\subset\R^3\times\R^3$ is measurable 
and of finite measure, then 
\[ 0\le\int\int_U f_n\,dx\,dv\to\int\int_U f_0\,dx\,dv,
   \quad n\to\infty, \]  
which yields $f_0\ge 0$ a.e. Similarly, 
\[ 0\le\int\int_U (1-f_n)\,dx\,dv\to\int\int_U (1-f_0)\,dx\,dv,
   \quad n\to\infty, \] 
for each such $U$ implies that $f_0\le 1$ a.e. 
Let $A\in {\rm SO}(3)$ and $X, V\subset\R^3$ be measurable and of finite measure. 
We have 
\begin{eqnarray*} 
   \int_X\int_V f_0(x, v)\,dx\,dv 
   & = & \lim_{n\to\infty}\int_X\int_V f_n(x, v)\,dx\,dv
   =\lim_{n\to\infty}\int_X\int_V f_n(Ax, Av)\,dx\,dv
   \\ & = & \lim_{n\to\infty}\int_{AX}\int_{AV} f_n(x, v)\,dx\,dv
   =\int_{AX}\int_{AV} f_0(x, v)\,dx\,dv
   \\ & = & \int_X\int_V f_0(Ax, Av)\,dx\,dv
\end{eqnarray*}  
as $n\to\infty$, and it follows that $f_0$ is spherically symmetric. 

If $P\ge P_0+1$, $R>0$ and $S>P+1$, then by (\ref{allPn}) for $n\in\N$ 
\[ \frac{2}{P^{1/4}}
   \ge\int_{\R^3}\int_{|v|\ge P+1}\sqrt{1+|v|^2}\,f_n\,dx\,dv
   \ge\int_{|x|\le R}\int_{P+1\le |v|\le S}\sqrt{1+|v|^2}\,f_n\,dx\,dv. \] 
Since ${\bf 1}_{|x|\le R}\sqrt{1+|v|^2}\,{\bf 1}_{P+1\le |v|\le S}
\in L^{1+1/k}(\R^3\times\R^3)^\ast$, we can pass to the limit $n\to\infty$ 
and get  
\[ \frac{2}{P^{1/4}}
   \ge\int_{|x|\le R}\int_{P+1\le |v|\le S}\sqrt{1+|v|^2}\,f_0\,dx\,dv \] 
for every $R>0$ and $S>P+1$. Taking the limits $R,S\to\infty$, 
we arrive at 
\begin{equation}\label{v0prop} 
   \int_{\R^3}\int_{|v|\ge P+1}\sqrt{1+|v|^2}\,f_0\,dx\,dv
   \le\frac{2}{P^{1/4}},\quad P\ge P_0+1.
\end{equation}  

Next fix $0\le r_1<r_2$. We claim that 
\begin{equation}\label{banst}
   \lim_{n\to\infty}\int_{r_1<|x|<r_2}\int_{\R^3}\sqrt{1+|v|^2}\,f_n\,dx\,dv
   =\int_{r_1<|x|<r_2}\int_{\R^3}\sqrt{1+|v|^2}\,f_0\,dx\,dv.
\end{equation}  
To establish (\ref{banst}), let $(P_N)$ be a sequence satisfying $P_N\ge P_0+1$ 
and $\lim_{N\to\infty} P_N=\infty$.
For $n, N\in\N$ we deduce from (\ref{allPn}) and (\ref{v0prop})   
\begin{eqnarray*} 
   \lefteqn{\bigg|\int_{r_1<|x|<r_2}\int_{\R^3}\sqrt{1+|v|^2} f_n\,dx\,dv
   -\int_{r_1<|x|<r_2}\int_{\R^3}\sqrt{1+|v|^2} f_0\,dx\,dv\bigg|}
   \\ & = & \bigg|\int_{r_1<|x|<r_2}\int_{|v|\le P_N+1}\sqrt{1+|v|^2} f_n\,dx\,dv 
   +\int_{r_1<|x|<r_2}\int_{|v|\ge P_N+1}\sqrt{1+|v|^2} f_n\,dx\,dv
   \\ & & \quad -\,\int_{r_1<|x|<r_2}\int_{|v|\le P_N+1}\sqrt{1+|v|^2} f_0\,dx\,dv
   -\int_{r_1<|x|<r_2}\int_{|v|\ge P_N+1}\sqrt{1+|v|^2} f_0\,dx\,dv\bigg| 
   \\ & \le & \bigg|\int_{r_1<|x|<r_2}\int_{|v|\le P_N+1}\sqrt{1+|v|^2} f_n\,dx\,dv 
   -\int_{r_1<|x|<r_2}\int_{|v|\le P_N+1}\sqrt{1+|v|^2} f_0\,dx\,dv\bigg|
   +\frac{4}{P_N^{1/4}}
   \\ & =: & \delta_{n, N}+\frac{4}{P_N^{1/4}}.     
\end{eqnarray*} 
Let $\eps>0$ be given. Then there is a $N_0\in\N$ such that 
$4P_{N_0}^{-1/4}\le\eps$. As $\lim_{n\to\infty}\delta_{n, N_0}=0$, 
we find an $n_0\in\N$ so that $\delta_{n, N_0}\le\eps$ for $n\ge n_0$, 
whence altogether we obtain $|\ldots|\le 2\eps$ for $n\ge n_0$, 
which completes the proof of (\ref{banst}). Taking $r_1=0$ and $r_2=r$, 
we deduce that in particular $\lim_{n\to\infty}m_{f_n}(r)=m_{f_0}(r)$ 
holds as $n\to\infty$, for $r\in [0, \infty[$. In addition, for $r>0$ and $\eps>0$ 
we get from (\ref{banst}) and $f_n\in {\cal A}$ 
\begin{eqnarray*} 
   4\pi\int_r^{r+\eps} s^2\rho_{f_0}(s)\,ds   
   & = & \int_{r<|x|<r+\eps}\int_{\R^3}\sqrt{1+|v|^2}\,f_0\,dx\,dv
   \\ & = & \lim_{n\to\infty}\int_{r<|x|<r+\eps}\int_{\R^3}\sqrt{1+|v|^2}\,f_n\,dx\,dv
   \\ & = & \lim_{n\to\infty}\int_{r<|x|<r+\eps}\rho_{f_n}(x)\,dx
   \\ & \le & \frac{4\pi\sigma_0}{3}\,((r+\eps)^3-r^3)
   =\frac{4\pi\sigma_0}{3}\,(3r^2\eps+3r\eps^2+\eps^3). 
\end{eqnarray*} 
Therefore, dividing by $\eps>0$ and taking the limit $\eps\to 0$, 
we see that for a.e.~$r>0$ we have $4\pi r^2\rho_{f_0}(r)\le 4\pi\sigma_0 r^2$, 
which means that $\rho_{f_0}(r)\le\sigma_0$. 

Lastly, to show that $\lambda_{f_n}(r)\to\lambda_{f_0}(r)$ pointwise, note first that
\begin{equation}\label{elamdiff} 
   e^{\lambda_{f_n}}-e^{\lambda_{f_0}}=\frac{2}{r}\,
   \frac{e^{2\lambda_{f_n}} e^{2\lambda_{f_0}}}{e^{\lambda_{f_n}}+e^{\lambda_{f_0}}}\,(m_{f_n}-m_{f_0}),
\end{equation} 
as is verified by a direct calculation. From above we know $\lim_{n\to\infty} m_{f_n}(r)=m_{f_0}(r)$ 
for every $r>0$. Since, owing to (\ref{cbetabd}), $1\le e^{2\lambda_{f_n}}\le c_\beta^2$ 
and $1\le e^{2\lambda_{f_0}}\le c_\beta^2$, the claim follows. 
{\hfill$\Box$}\bigskip

%%%%%%%%%%%%%%%%%%%%%%%%%%%%%%%%%%%%%%%%%%%%%%%%%%%%%%%%%%%%%%%%%%%%%%%%%%%%%%%%%%%%%%%%%

\section{The Casimir-binding-energy condition and existence of a minimizer}
\label{sec-spatial-part}
\setcounter{equation}{0}

We continue to consider
\begin{equation}\label{calH2} 
   {\cal D}(f)=\int_{\R^3}\int_{\R^3} e^{\lambda}\,(\chi(f)-f)\,dx\,dv,
\end{equation} 
which is (\ref{calH1}). 

First we study the scaling properties of the system; also see \cite{GRS}. 
Henceforth we will assume that $k\in ]0, 2]$. 

\begin{lemma}\label{scalem}  
Suppose that $f=f(x, v)\ge 0$ is measurable 
and spherically symmetric such that 
$\int_{\R^3}\int_{\R^3}\sqrt{1+|v|^2} f\,dx\,dv=M$. For $\gamma>0$ define 
\[ f_\gamma(x, v)=\gamma^2 f(\gamma x, v). \] 
Then 
\[ M_\gamma=\gamma^{-1} M,\quad\rho_\gamma(r)=\gamma^2\rho(\gamma r),
   \quad m_\gamma(r)=\gamma^{-1}m(\gamma r),\quad   
   \lambda_\gamma(r)=\lambda(\gamma r), \] 
and if furthermore $\gamma\le 1$, then also 
\[ {\cal D}(f_\gamma)\le\gamma^{-1}\,{\cal D}(f). \] 
\end{lemma} 
{\bf Proof\,:} By changing variables as $\tilde{x}=\gamma x$, $d\tilde{x}=\gamma^3 dx$, 
we calculate 
\[ M_\gamma=\int_{\R^3}\int_{\R^3}\sqrt{1+|v|^2} f_\gamma\,dx\,dv 
   =\gamma^2\gamma^{-3}\int_{\R^3}\int_{\R^3}\sqrt{1+|v|^2} f(\tilde{x}, v)\,d\tilde{x}\,dv
   =\gamma^{-1} M. \] 
Moreover, 
\[ \rho_\gamma(r)=\int_{\R^3}\sqrt{1+|v|^2} f_\gamma(x, v)\,dv=\gamma^2\rho(\gamma r)  \] 
yields 
\[ m_\gamma(r)=4\pi\int_0^r s^2\rho_\gamma(s)\,ds
   =4\pi\gamma^2\int_0^r s^2\rho(\gamma s)\,ds
   =4\pi\gamma^{-1}\int_0^{\gamma r} \tau^2\rho(\tau)\,d\tau
   =\gamma^{-1}m(\gamma r). \] 
This in turn leads to 
\[ e^{-2\lambda_\gamma(r)}=1-\frac{2m_\gamma(r)}{r}
   =1-\frac{m(\gamma r)}{\gamma r}=e^{-2\lambda(\gamma r)}, \] 
so that $\lambda_\gamma(r)=\lambda(\gamma r)$. 
It follows that 
\begin{eqnarray*} 
   {\cal D}(f_\gamma) & = & \int_{\R^3}\int_{\R^3} e^{\lambda_\gamma}\,(\chi(f_\gamma)-f_\gamma)\,dx\,dv 
   \\ & = & \int_{\R^3}\int_{\R^3} e^{\lambda(\gamma r)}
   \,\Big(\frac{k}{k+1}\,\gamma^{2(1+\frac{1}{k})}\,f(\gamma x, v)^{1+\frac{1}{k}}
   -\gamma^2 f(\gamma x, v)\Big)\,dx\,dv
   \\ & = & \int_{\R^3}\int_{\R^3} e^{\lambda(\tilde{r})}
   \,\Big(\frac{k}{k+1}\,\gamma^{\frac{2}{k}-1}\,f(\tilde{x}, v)^{1+\frac{1}{k}}
   -\gamma^{-1} f(\tilde{x}, v)\Big)\,d\tilde{x}\,dv
   \\ & = & {\cal D}(f)+\int_{\R^3}\int_{\R^3} e^{\lambda(\tilde{r})}
   \,\Big(\frac{k}{k+1}\,[\gamma^{\frac{2}{k}-1}-1]\,f(\tilde{x}, v)^{1+\frac{1}{k}}
   -[\gamma^{-1}-1] f(\tilde{x}, v)\Big)\,d\tilde{x}\,dv. 
\end{eqnarray*} 
Since $k\in ]0, 2]$ we have $\frac{2}{k}-1\ge 0$, and hence $\gamma^{\frac{2}{k}-1}-1\le 0$ 
due to $\gamma\le 1$. Therefore we obtain 
\begin{eqnarray*} 
   {\cal D}(f_\gamma) & \le & {\cal D}(f)+(\gamma^{-1}-1)\int_{\R^3}\int_{\R^3} e^{\lambda}\,(-f)\,dx\,dv
   \\ & \le & {\cal D}(f)+(\gamma^{-1}-1)\int_{\R^3}\int_{\R^3} e^{\lambda}\,(\chi(f)-f)\,dx\,dv
   \\ & = & \gamma^{-1}\,{\cal D}(f),  
\end{eqnarray*} 
as claimed. {\hfill$\Box$}\bigskip

\begin{lemma}\label{topp20} 
Let $f\in {\cal A}$ and $\eps, \eta, R>0$ be such that $\frac{2M}{R}\le 1-\frac{1}{(1+\eta)^2}$, 
\[ f\le 1,\quad {\cal D}(f)\le I+\eps\quad\mbox{and}
   \quad\Delta M=\int_{|x|\ge R}\int_{\R^3}\sqrt{1+|v|^2}\,f\,dx\,dv\in ]0, M[ \] 
are satisfied. If $f_1={\bf 1}_{\{|x|\le R\}} f$, then 
\[ \frac{\Delta M}{M-\Delta M}\,|{\cal D}(f_1)|
   \le (1+\eta)\Delta M+\eps. \] 
\end{lemma} 
{\bf Proof\,:} Let $f_2={\bf 1}_{\{|x|\ge R\}} f$. 
Since $f\le 1$ we have $\chi(f)-f\le 0$ pointwise, 
recall (\ref{sith}). Therefore in particular $\chi(f_1)-f_1\le 0$ pointwise 
and $\chi(f_2)-f_2\le 0$ pointwise. If $r\le R$, then 
\begin{eqnarray}\label{m1m}  
   m_{f_1}(r) & = & 4\pi\int_0^r s^2\rho_{f_1}(s)\,ds
   =\int_{|x|\le r}\int_{\R^3}\sqrt{1+|v|^2}\,f_1(x, v)\,dx\,dv
   \nonumber
   \\ & = & \int_{|x|\le r}\int_{\R^3}\sqrt{1+|v|^2}\,f(x, v)\,dx\,dv=m_f(r), 
\end{eqnarray} 
and hence also 
\[ e^{\lambda_{f_1}(r)}=e^{\lambda_f(r)},\quad r\le R. \] 
This gives 
\begin{eqnarray*} 
   |{\cal D}(f)| & = & \int_{\R^3}\int_{\R^3} e^{\lambda_f(r)}\,(f-\chi(f))\,dx\,dv 
   \\ & = & \int_{|x|\le R}\int_{\R^3} e^{\lambda_f(r)}\,(f-\chi(f))\,dx\,dv 
   +\int_{|x|\ge R}\int_{\R^3} e^{\lambda_f(r)}\,(f-\chi(f))\,dx\,dv 
   \\ & = & \int_{|x|\le R}\int_{\R^3} e^{\lambda_{f_1}(r)}\,(f_1-\chi(f_1))\,dx\,dv 
   +\int_{|x|\ge R}\int_{\R^3} e^{\lambda_f(r)}\,(f-\chi(f))\,dx\,dv 
   \\ & = & |{\cal D}(f_1)|+\int_{|x|\ge R}\int_{\R^3} e^{\lambda_f(r)}\,(f-\chi(f))\,dx\,dv. 
\end{eqnarray*} 
If $r\ge R$, then $\frac{2M}{R}\le 1-\frac{1}{(1+\eta)^2}$ leads to    
\[ e^{\lambda_f(r)}=\frac{1}{\sqrt{1-\frac{2m_f(r)}{r}}}
   \le\frac{1}{\sqrt{1-\frac{2M}{R}}}\le 1+\eta. \] 
It follows that 
\begin{eqnarray*} 
   \int_{|x|\ge R}\int_{\R^3} e^{\lambda_f(r)}\,(f-\chi(f))\,dx\,dv
   & \le & (1+\eta)\int_{|x|\ge R}\int_{\R^3} f\,dx\,dv
   \\ & \le & (1+\eta)\int_{|x|\ge R}\int_{\R^3}\sqrt{1+|v|^2}\,f\,dx\,dv
   \le (1+\eta)\Delta M.
\end{eqnarray*} 
As a consequence, 
\begin{equation}\label{aac} 
   |{\cal D}(f)|\le |{\cal D}(f_1)|+(1+\eta)\Delta M.
\end{equation} 
Next we apply the scaling property from Lemma \ref{scalem}
to $f_1$, $M$ replaced by the total mass $M-\Delta M$ of $f_1$ 
and $\gamma=\frac{M-\Delta M}{M}\in ]0, 1[$ to obtain 
a spherically symmetric function $f_{1, \gamma}\ge 0$ such that  
\[ \int_{\R^3}\int_{\R^3}\sqrt{1+|v|^2} f_{1,\gamma}\,dx\,dv=\gamma^{-1} (M-\Delta M)=M \] 
and also
\[ {\cal D}(f_{1,\gamma})\le\frac{M}{M-\Delta M}\,{\cal D}(f_1)\le 0. \] 
From $\rho_{f_{1,\gamma}}(r)=\gamma^2\rho_{f_1}(\gamma r)\le\rho_{f_1}(\gamma r)\le\sigma_0$ 
we deduce that $f_{1,\gamma}\in {\cal A}$, so that 
\[ {\cal D}(f)-\eps\le I\le {\cal D}(f_{1,\gamma})
   \le\frac{M}{M-\Delta M}\,{\cal D}(f_1). \]  
Taking minus signs, this can be rewritten as 
\[ \frac{M}{M-\Delta M}\,|{\cal D}(f_1)|\le |{\cal D}(f)|+\eps. \] 
Using (\ref{aac}), we arrive at 
\[ \frac{\Delta M}{M-\Delta M}\,|{\cal D}(f_1)|
   \le (1+\eta)\Delta M+\eps, \] 
as was to be shown. 
{\hfill$\Box$}\bigskip

Now we turn to the Casimir-binding-energy condition (CBEC), 
as stated in Definition \ref{sbec}. 

\begin{remark}\label{sbech}{\rm (a) Observe that (CBEC) holds if and only if 
there is a function 
$h\in {\cal A}$ such that ${\cal D}(h)<0$ and 
\[ |{\cal D}(h)|>M. \]    
\smallskip 

\noindent 
(b) According to (a) and Lemma \ref{Ikle}, (CBEC) is verified 
in particular for $k\in ]0, 2]$. 
}
\end{remark} 

\begin{cor}\label{SBECcor} 
Let $h$ be as in Remark \ref{sbech}(a) and $\delta>0$ 
such that $|{\cal D}(h)|\ge M+\delta$. Next define 
\[ \eta=\min\Big\{1, \frac{\delta}{6M}\Big\} \] 
and take 
\begin{equation}\label{Rcond} 
   R\ge\max\Big\{\frac{8M}{3}, \frac{2M}{1-\frac{1}{(1+\frac{\delta}{6M})^2}}\Big\}.
\end{equation}  
If $g\in {\cal A}$ satisfies 
\[ g\le 1,\quad {\cal D}(g)\le I+\eta (\Delta M)
   \quad\mbox{and}\quad\Delta M\in ]0, M[, \] 
where 
\[ \Delta M=\int_{|x|\ge R}\int_{\R^3}\sqrt{1+|v|^2}\,g\,dx\,dv, \] 
then 
\[ \Delta M\ge\frac{\delta}{18}. \] 
In particular, if $(f_n)\subset {\cal A}$ is a minimizing sequence 
such that $f_n\le 1$ and
\[ \lim_{n\to\infty} {(\Delta M)}_n=\eps_0 \] 
for some $\eps_0\in ]0, M[$, where 
\begin{equation}\label{Rncond} 
   {(\Delta M)}_n=\int_{|x|\ge R_n}\int_{\R^3}\sqrt{1+|v|^2}\,f_n\,dx\,dv
   \quad\mbox{with}\quad R_n\ge\max\Big\{\frac{8M}{3}, \frac{2M}{1-\frac{1}{(1+\frac{\delta}{6M})^2}}\Big\},
\end{equation}  
then 
\[ \eps_0\ge\frac{\delta}{18}. \] 
\end{cor} 
{\bf Proof\,:} We can apply Lemma \ref{topp20} with $f=g$ and $\eps=\eta (\Delta M)$, 
since $\frac{2M}{R}\le 1-\frac{1}{(1+\eta)^2}$. It follows that 
\[ \frac{\Delta M}{M-\Delta M}\,|{\cal D}(g_1)|
   \le (1+\eta)\Delta M+\eps=(1+2\eta)\Delta M \] 
for $g_1={\bf 1}_{\{|x|\le R\}}\,g$, and hence  
\[ |{\cal D}(g_1)|\le (M-\Delta M)(1+2\eta). \] 
Now we distinguish two cases. 
Case 1: $|{\cal D}(h)|\le |{\cal D}(g_1)|+\frac{\delta}{6}$. Then 
\begin{eqnarray*} 
   M+\delta & \le & |\mathcal{D}(h)|\le |{\cal D}(g_1)|+\frac{\delta}{6} 
   \le (M-\Delta M)(1+2\eta)+\frac{\delta}{6}
   \\ & \le & M(1+2\eta)+\frac{\delta}{6}\le M+\frac{\delta}{2}, 
\end{eqnarray*} 
which is a contradiction. 
Case 2: $|{\cal D}(h)|\ge |{\cal D}(g_1)|+\frac{\delta}{6}$. 
If $r\le R$, then 
\[ m_{g_1}(r)=4\pi\int_0^r s^2\rho_{g_1}(s)\,ds
   =\int_{|x|\le r}\int_{\R^3} g_1(x, v)\,dx\,dv
   =\int_{|x|\le r}\int_{\R^3} g(x, v)\,dx\,dv
   =m_g(r), \] 
and hence also $\lambda_{g_1}(r)=\lambda_g(r)$ for $r\le R$.  
Therefore we obtain 
\begin{eqnarray*} 
   {\cal D}(h) & \le & {\cal D}(g_1)-\frac{\delta}{6}
   \\ & = & \int_{|x|\le R}\int_{\R^3} e^{\lambda_{g_1}} (\chi(g_1)-g_1)\,dx\,dv-\frac{\delta}{6}
   \\ & = & \int_{|x|\le R}\int_{\R^3} e^{\lambda_g} (\chi(g)-g)\,dx\,dv-\frac{\delta}{6}
   \\ & = & {\cal D}(g)+\int_{|x|\ge R}\int_{\R^3} e^{\lambda_g} (g-\chi(g))\,dx\,dv-\frac{\delta}{6} 
   \\ & \le & {\cal D}(g)+\frac{1}{\sqrt{1-\frac{2M}{R}}}
   \int_{|x|\ge R}\int_{\R^3} g\,dx\,dv-\frac{\delta}{6} 
   \\ & \le & {\cal D}(g)+(1+\eta)\Delta M-\frac{\delta}{6} 
   \\ & \le & I+\eta (\Delta M)+(1+\eta)\Delta M-\frac{\delta}{6}
   \\ & = & I+(1+2\eta)\Delta M-\frac{\delta}{6}
   \\ & \le & I+3\Delta M-\frac{\delta}{6}. 
\end{eqnarray*}  
Since $h\in {\cal A}$, we have $I\le {\cal D}(h)$, which implies that 
$\Delta M\ge\frac{\delta}{18}$. 

For the second part, $R_n$ satisfies (\ref{Rcond}) by hypothesis. 
Due to the fact that $(f_n)$ is a minimizing sequence, 
we have $\lim_{n\to\infty} {\cal D}(f_n)=I$. 
Then ${(\Delta M)}_n\to\eps_0>0$ shows that ${\cal D}(f_n)\le I+\eta{(\Delta M)}_n$ 
is verified for $n$ large enough. In addition, $\eps_0\in ]0, M[$ 
yields that also ${(\Delta M)}_n\in ]0, M[$ for $n$ sufficiently large.  	
Thus the first part applies to deduce that ${(\Delta M)}_n\ge\frac{\delta}{18}$ 
for large $n$, so that $\eps_0\ge\frac{\delta}{18}$ is found in the limit $n\to\infty$. 
{\hfill$\Box$}\bigskip

The last part of the previous result can be sharpened as follows. 

\begin{cor}\label{SBECreal} 
Let $h$ be as in Remark \ref{sbech}(a) and $\delta>0$ 
such that $|{\cal D}(h)|\ge M+\delta$. Let $(R_n)$ be a sequence such that 
\begin{equation}\label{Rncond2} 
   R_n\ge\max\Big\{\frac{8M}{3}, \frac{2M}{1-\frac{1}{(1+\frac{\delta}{6M})^2}}\Big\},
   \quad n\in\N.  
\end{equation} 
Then there is no minimizing sequence $(f_n)\subset {\cal A}$ such that $f_n\le 1$ and
\[ \lim_{n\to\infty} {(\Delta M)}_n\in ]0, M[, \] 
where 
\[ {(\Delta M)}_n=\int_{|x|\ge R_n}\int_{\R^3}\sqrt{1+|v|^2}\,f_n\,dx\,dv. \] 
\end{cor} 
{\bf Proof\,:} Since $\delta>0$ has the property that $|{\cal D}(h)|\ge M+\delta$, 
there is no loss of generality in assuming that $\delta<20 M$, 
as otherwise $\delta$ could be replaced by a smaller number. 
Assuming the contrary of what is stated, 
let $\eps_0=\lim_{n\to\infty} {(\Delta M)}_n\in ]0, M[$. 
We distinguish two cases. Case 1: $\eps_0\le\frac{\delta}{20}$. 
From Corollary \ref{SBECcor} we however deduce that $\frac{\delta}{18}\le\eps_0\le\frac{\delta}{20}$, 
which is a contradiction. Case 2: $\eps_0>\frac{\delta}{20}$. Then we may assume that 
\[ {(\Delta M)}_n=\int_{|x|\ge R_n}\int_{\R^3}\sqrt{1+|v|^2}\,f_n\,dx\,dv>\frac{\delta}{20}, 
   \quad n\in\N. \] 
For every fixed $n\in\N$ we consider the continuous function $\kappa_n(R)
=\int_{|x|\ge R}\int_{\R^3}\sqrt{1+|v|^2}\,f_n\,dx\,dv$. 
Then $\kappa_n(R_n)>\frac{\delta}{20}$ and $\lim_{R\to\infty}\kappa_n(R)=0$. 
Hence there is a $\tilde{R}_n>R_n$ such that $\kappa_n(\tilde{R}_n)=\frac{\delta}{20}$. For 
\[ {(\widetilde{\Delta M})}_n=\int_{|x|\ge\tilde{R}_n}\int_{\R^3}\sqrt{1+|v|^2}\,f_n\,dx\,dv \]  
we then get ${(\widetilde{\Delta M})}_n=\kappa_n(\tilde{R}_n)=\frac{\delta}{20}$ for $n\in\N$. 
Furthermore, since $\tilde{R}_n>R_n$, also $\tilde{R}_n$ satisfies (\ref{Rncond}). 
As $\frac{\delta}{20}<M$ by the remark at the beginning of the argument, 
we can once again apply Corollary \ref{SBECcor} to yield $\frac{\delta}{20}\ge\frac{\delta}{18}$, 
which is impossible. 
{\hfill$\Box$}\bigskip

From now on we take $(f_n)\subset {\cal A}$ to be a minimizing sequence for ${\cal D}$ 
with the additional property (\ref{allPn}), and $0\le f_n\le 1$ for $n\in\N$, 
that led to the weak limit $f_0$; see the paragraphs before 
and following Lemma \ref{1+epsbd}. 

\begin{lemma} For every $R>0$ we have 
\begin{equation}\label{eram} 
   \lim_{n\to\infty}\int_{|x|\le R}\int_{\R^3}\sqrt{1+|v|^2}\,f_n\,dx\,dv
   =\int_{|x|\le R}\int_{\R^3}\sqrt{1+|v|^2}\,f_0\,dx\,dv
\end{equation} 
and 
\begin{equation}\label{eram2} 
   \lim_{n\to\infty}\int_{|x|\le R}\int_{\R^3}\varphi\,f_n\,dx\,dv
   =\int_{|x|\le R}\int_{\R^3}\varphi\,f_0\,dx\,dv
\end{equation} 
for every bounded and continuous function $\varphi=\varphi(x, v)\ge 0$. 
\end{lemma} 
{\bf Proof\,:} For $n\in\N$ and $P\ge P_0+1$ we use (\ref{allPn}) and Lemma \ref{f0props} 
to bound 
\begin{eqnarray*} 
   \lefteqn{ \bigg|\int_{|x|\le R}\int_{\R^3}\sqrt{1+|v|^2}\,f_n\,dx\,dv
   -\int_{|x|\le R}\int_{\R^3}\sqrt{1+|v|^2}\,f_0\,dx\,dv\bigg|} 
   \\ & = & 
   \bigg|\int_{|x|\le R}\int_{|v|\le P+1}\sqrt{1+|v|^2}\,f_n\,dx\,dv
   +\int_{|x|\le R}\int_{|v|\ge P+1}\sqrt{1+|v|^2}\,f_n\,dx\,dv
   \\ & & -\,\int_{|x|\le R}\int_{|v|\le P+1}\sqrt{1+|v|^2}\,f_0\,dx\,dv
   -\int_{|x|\le R}\int_{|v|\ge P+1}\sqrt{1+|v|^2}\,f_0\,dx\,dv\bigg|
   \\ & \le & \bigg|\int_{|x|\le R}\int_{|v|\le P+1}\sqrt{1+|v|^2}\,f_n\,dx\,dv
   -\int_{|x|\le R}\int_{|v|\le P+1}\sqrt{1+|v|^2}\,f_0\,dx\,dv\bigg|
   +\frac{4}{P^{1/4}}. 
\end{eqnarray*} 
Given $\eps>0$ choose $P\ge P_0+1$ large enough so that $4P^{-{1/4}}\le\eps/2$. 
Since $|\ldots|\to 0$ as $n\to\infty$ for the first term on the right-hand side 
due to the weak convergence, we find $n_0\in\N$ with the property that 
$|\ldots|\le\eps/2$ for $n\ge n_0$, which establishes (\ref{eram}). 

To show (\ref{eram2}), one can proceed in an analogous way, using that this time 
\begin{eqnarray*} 
   \lefteqn{\bigg|\int_{|x|\le R}\int_{|v|\ge P+1}\varphi\,f_n\,dx\,dv
   -\int_{|x|\le R}\int_{|v|\ge P+1}\varphi\,f_0\,dx\,dv\bigg|} 
   \\ & \le & \int_{|x|\le R}\int_{|v|\ge P+1}\varphi\,f_n\,dx\,dv
   +\int_{|x|\le R}\int_{|v|\ge P+1}\varphi\,f_0\,dx\,dv
   \\ & \le & {\|\varphi\|}_\infty\int_{|x|\le R}\int_{|v|\ge P+1}\sqrt{1+|v|^2}\,f_n\,dx\,dv
   +{\|\varphi\|}_\infty\int_{|x|\le R}\int_{|v|\ge P+1}\sqrt{1+|v|^2}\,f_0\,dx\,dv
   \\ & \le & \frac{4{\|\varphi\|}_\infty}{P^{1/4}}
\end{eqnarray*} 
for the error term. For the other term, the weak convergence may still be used, 
owing to the fact that ${\bf 1}_{\{|x|\le R\}}{\bf 1}_{\{|v|\le P+1\}}\varphi
\in {L^{1+1/k}(\R^3\times\R^3)}^\ast$. 
{\hfill$\Box$}\bigskip 

Now we are going to show that (CBEC) implies that $f_0$ has compact support in $x$. 

\begin{lemma}\label{f0cpt} 
Let $f_0\in L^{1+1/k}(\R^3\times\R^3)$ 
be the limiting function as obtained in Lemma \ref{f0props}. 
Then there exists a $R_0>0$ such that 
\[ \int_{\R^3}\int_{|x|\ge R_0}\sqrt{1+|v|^2} f_0\,dx\,dv=0, \] 
which means that $f_0(x, v)=0$ for a.e.~$|x|\ge R_0$ and $v\in\R^3$. 
\end{lemma} 
{\bf Proof\,:} Let $\kappa(R)=\int_{|x|\ge R}\int_{\R^3}\sqrt{1+|v|^2}\,f_0\,dx\,dv$ 
and assume that $\kappa(R)>0$ for all $R>0$. Let $h$ be as in Remark \ref{sbech}(a) and $\delta>0$ 
such that $|{\cal D}(h)|\ge M+\delta$. Since $\lim_{R\to\infty}\kappa(R)=0$, there is 
\[ R_0\ge\max\Big\{\frac{8M}{3}, \frac{2M}{1-\frac{1}{(1+\frac{\delta}{6M})^2}}\Big\} \] 
so that 
\[ 0<\eps_0:=\kappa(R_0)\le\frac{M}{2}. \] 
In particular, $\eps_0\in ]0, M[$, and in the present case, the relation (\ref{eram}) can be rewritten as 
\[ \lim_{n\to\infty}\int_{|x|\le R_0}\int_{\R^3}\sqrt{1+|v|^2}\,f_n\,dx\,dv=M-\eps_0. \]  
As $f_n$ has mass $M$, this in turn implies 
\[ \lim_{n\to\infty}\int_{|x|\ge R_0}\int_{\R^3}\sqrt{1+|v|^2}\,f_n\,dx\,dv=\eps_0. \] 
Since $(f_n)$ is a minimizing sequence, we obtain a contradiction to Corollary \ref{SBECreal}, 
upon taking $R_n=R_0$.  
{\hfill$\Box$}\bigskip

\begin{lemma}\label{f0incalA} 
Let $f_0\in L^{1+1/k}(\R^3\times\R^3)$ 
be the limiting function as obtained in Lemma \ref{f0props}. Then $f_0\neq 0$ 
and $f_0\in {\cal A}$. 
\end{lemma} 
{\bf Proof\,:} First we are going to show that $f_0\neq 0$. Otherwise we would get 
\[ \lim_{n\to\infty}\int_{|x|\le R}\int_{\R^3}\sqrt{1+|v|^2}\,f_n\,dx\,dv=0 \] 
for every $R>0$ from (\ref{eram}). Let $h$ be as in Remark \ref{sbech}(a) and $\delta>0$ 
such that $|{\cal D}(h)|\ge M+\delta$. In particular, $h\in {\cal A}$ implies that 
$I\le {\cal D}(h)\le -M-\delta$. Fix 
\[ R=\frac{2M}{1-\frac{1}{(1+\frac{\delta}{2M})^2}}>2M. \] 
If $n\in\N$, then by (\ref{cbetabd}) 
\begin{eqnarray*} 
   |{\cal D}(f_n)| & = & 
   \int_{\R^3}\int_{\R^3} e^{\lambda_n}\,(f_n-\chi(f_n))\,dx\,dv 
   \le\int_{\R^3}\int_{\R^3} e^{\lambda_n}\,f_n\,dx\,dv
   \\ & \le & c_\beta\int_{|x|\le R}\int_{\R^3}\sqrt{1+|v|^2}\,f_n\,dx\,dv 
   +\frac{1}{\sqrt{1-\frac{2M}{R}}}\int_{|x|\ge R}\int_{\R^3} f_n\,dx\,dv
   \\ & \le & c_\beta\int_{|x|\le R}\int_{\R^3}\sqrt{1+|v|^2}\,f_n\,dx\,dv 
   +\frac{1}{\sqrt{1-\frac{2M}{R}}}\,M
   \\ & = & c_\beta\int_{|x|\le R}\int_{\R^3}\sqrt{1+|v|^2}\,f_n\,dx\,dv 
   +\Big(1+\frac{\delta}{2M}\Big)M. 
\end{eqnarray*} 
Passing to the limit $n\to\infty$, we deduce 
\[ M+\delta\le |I|\le M+\frac{\delta}{2},  \] 
which is a contradiction. 
\smallskip 

\noindent 
In order to establish that $f_0\in {\cal A}$, due to Lemma \ref{f0props} 
the only thing that remains to be verified is that $f_0$ has mass $M$. 
Suppose that $m=\int_{\R^3}\int_{\R^3}\sqrt{1+|v|^2} f_0\,dx\,dv<M$. 
Due to $f_0\neq 0$, then in fact $m\in ]0, M[$. Let $R_0>0$ be such that 
$f_0(x, v)=0$ for a.e.~$|x|\ge R_0$ and $v\in\R^3$; this is in accordance with Lemma \ref{f0cpt}. 
By (\ref{eram}) it follows that for every $R\ge R_0$ we have 
\begin{eqnarray*} 
   \lim_{n\to\infty}\int_{|x|\le R}\int_{\R^3}\sqrt{1+|v|^2}\,f_n\,dx\,dv
   & = & \int_{|x|\le R}\int_{\R^3}\sqrt{1+|v|^2}\,f_0\,dx\,dv
   \\ & = & \int_{|x|\le R_0}\int_{\R^3}\sqrt{1+|v|^2}\,f_0\,dx\,dv=m.
\end{eqnarray*}  
Since every $f_n$ has mass $M$, we deduce that, for every $R\ge R_0$,  
\[ \lim_{n\to\infty}\int_{|x|\ge R}\int_{\R^3}\sqrt{1+|v|^2}\,f_n\,dx\,dv
   =\eps_0 \] 
for $\eps_0=M-m\in ]0, M[$. In particular, if $\tilde{R}_0\ge R_0$ satisfies 
\[ \tilde{R}_0\ge\max\Big\{\frac{8M}{3}, \frac{2M}{1-\frac{1}{(1+\frac{\delta}{6M})^2}}\Big\}, \]
then 
\[ \lim_{n\to\infty}\int_{|x|\ge\tilde{R}_0}\int_{\R^3}\sqrt{1+|v|^2}\,f_n\,dx\,dv
   =\eps_0. \]
Since $(f_n)$ is a minimizing sequence, we may now apply Corollary \ref{SBECreal} 
with $R_n=\tilde{R}_0$ to get a contradiction. 
{\hfill$\Box$}\bigskip

\begin{lemma}\label{2ndpartconv}
Let $f_0\in L^{1+1/k}(\R^3\times\R^3)$ 
be the limiting function as obtained in Lemma \ref{f0props}. Then 
\[ \lim_{n\to\infty}\int_{\R^3}\int_{\R^3} e^{\lambda_{f_n}} f_n\,dx\,dv
   =\int_{\R^3}\int_{\R^3} e^{\lambda_{f_0}} f_0\,dx\,dv, \]
up to taking a further subsequence.   
\end{lemma} 
{\bf Proof\,:} First we will show that 
\begin{equation}\label{gaspec} 
   \lim_{n\to\infty} \int_{\R^3}\int_{\R^3} |e^{\lambda_{f_n}}-e^{\lambda_{f_0}}|\,f_n\,dx\,dv=0.
\end{equation}  
For, let $R>0$ and $P\ge P_0+1$. 
Then due to $f_n\le 1$, by (\ref{elamdiff}), (\ref{cbetabd}) 
and (\ref{allPn}) [remember that $\tilde{f}_n$ had been renamed to $f_n$], 
\begin{eqnarray}\label{dihu} 
   \lefteqn{\int_{\R^3}\int_{\R^3} |e^{\lambda_{f_n}}-e^{\lambda_{f_0}}|\,f_n\,dx\,dv}
   \nonumber
   \\ & = & \int_{|x|\le R}\int_{|v|\le P+1} |e^{\lambda_{f_n}}-e^{\lambda_{f_0}}|\,f_n\,dx\,dv
   +\int_{|x|\ge R}\int_{|v|\le P+1} |e^{\lambda_{f_n}}-e^{\lambda_{f_0}}|\,f_n\,dx\,dv
   \nonumber
   \\ & & +\,\int_{\R^3}\int_{|v|\ge P+1} |e^{\lambda_{f_n}}-e^{\lambda_{f_0}}|\,f_n\,dx\,dv
   \nonumber
   \\ & \le & \int_{|x|\le R}\int_{|v|\le P+1} |e^{\lambda_{f_n}}-e^{\lambda_{f_0}}|\,dx\,dv
   +\int_{|x|\ge R}\int_{|v|\le P+1}\frac{2}{r}\,
   \frac{e^{2\lambda_{f_n}} e^{2\lambda_{f_0}}}{e^{\lambda_{f_n}}+e^{\lambda_{f_0}}}
   \,|m_{f_n}-m_{f_0}|\,f_n\,dx\,dv
   \nonumber
   \\ & & +\,2c_\beta\int_{\R^3}\int_{|v|\ge P+1} \sqrt{1+|v|^2}\,f_n\,dx\,dv
   \nonumber
   \\ & \le & \int_{|x|\le R}\int_{|v|\le P+1} |e^{\lambda_{f_n}}-e^{\lambda_{f_0}}|\,dx\,dv
   +\frac{2M}{R}\,c_\beta^4\int_{|x|\ge R}\int_{|v|\le P+1} f_n\,dx\,dv
   +\frac{4c_\beta}{P^{1/4}}
   \nonumber
   \\ & \le & \int_{|x|\le R}\int_{|v|\le P+1} |e^{\lambda_{f_n}}-e^{\lambda_{f_0}}|\,dx\,dv
   +\frac{2M^2}{R}\,c_\beta^4+\frac{4c_\beta}{P^{1/4}}. 
\end{eqnarray} 
Given $\eps>0$ we fix $R>0$ and $P\ge P_0+1$ so large that 
the last two terms add up to something less than $\eps/2$. 
The dominated convergence theorem implies that 
\begin{equation}\label{doco} 
   \lim_{n\to\infty}\int_{|x|\le R}\int_{|v|\le P+1} |e^{\lambda_{f_n}}-e^{\lambda_{f_0}}|\,dx\,dv=0.
\end{equation}  
In fact, we have $e^{\lambda_{f_n}(r)}\to e^{\lambda_{f_0}(r)}$ as $n\to\infty$ for $r\in [0, \infty[$ 
pointwise by Lemma \ref{f0props} and moreover $|e^{\lambda_{f_n}}-e^{\lambda_{f_0}}|
\le 2c_\beta$ by (\ref{cbetabd}). Thus (\ref{doco}) holds, and accordingly we may choose 
$n_0\in\N$ such that the first term in (\ref{dihu}) is less that $\eps/2$ for every $n\ge n_0$. 
This completes the argument for verifying (\ref{gaspec}). 

Therefore in order to prove the lemma, it remains to show that 
\begin{equation}\label{ropla} 
   \lim_{n\to\infty} \int_{\R^3}\int_{\R^3} e^{\lambda_{f_0}} f_n\,dx\,dv
   =\int_{\R^3}\int_{\R^3} e^{\lambda_{f_0}} f_0\,dx\,dv.
\end{equation}   
From (\ref{eram2}) and (\ref{cbetabd}) we deduce that 
\[ \lim_{n\to\infty}\int_{|x|\le R}\int_{\R^3} e^{\lambda_{f_0}} f_n\,dx\,dv
   =\int_{|x|\le R}\int_{\R^3} e^{\lambda_{f_0}} f_0\,dx\,dv \]  
holds for every $R>0$. Choose $R_0>0$ such that 
$f_0(x, v)=0$ for a.e.~$|x|\ge R_0$ and $v\in\R^3$; this is possible owing to Lemma \ref{f0cpt}. 
If $R\ge R_0$, then 
\[ \lim_{n\to\infty}\int_{|x|\le R}\int_{\R^3} e^{\lambda_{f_0}} f_n\,dx\,dv
   =\int_{|x|\le R}\int_{\R^3} e^{\lambda_{f_0}} f_0\,dx\,dv
   =\int_{\R^3}\int_{\R^3} e^{\lambda_{f_0}} f_0\,dx\,dv=:\theta_0, \] 
and $\theta_0>0$. Let $(\eps_m)$ be a sequence such that $\eps_m\downarrow 0$ 
as $m\to\infty$. Hence for every $m\in\N$ and $R\ge R_0$ there is $n(m, R)\in\N$ such that 
\[ \bigg|\int_{|x|\le R}\int_{\R^3} e^{\lambda_{f_0}} f_n\,dx\,dv-\theta_0\bigg|\le\eps_m,
   \quad n\ge n(m, R). \] 
Taking $R=m$ (for $m$ sufficiently large, say $m\ge m_0$) and defining $n_m=n(m, m)+m$, 
we then obtain 
\[ \bigg|\int_{|x|\le m}\int_{\R^3} e^{\lambda_{f_0}} f_{n_m}\,dx\,dv-\theta_0\bigg|\le\eps_m,
   \quad m\ge m_0, \] 
as well as $n_m\to\infty$ as $m\to\infty$. 
Now we need to distinguish two cases.  
Case 1: There is a further subsequence $(f_{n_{m_k}})$ of $(f_{n_m})$ such that  
\[ \lim_{k\to\infty}\int_{|x|\ge m_k}\int_{\R^3} e^{\lambda_{f_0}} f_{n_{m_k}}\,dx\,dv=0. \]  
Then 
\begin{eqnarray*} 
   \lefteqn{\bigg|\int_{\R^3}\int_{\R^3} e^{\lambda_{f_0}} f_{n_{m_k}}\,dx\,dv
   -\int_{\R^3}\int_{\R^3} e^{\lambda_{f_0}} f_0\,dx\,dv\bigg|} 
   \\ & \le & \int_{|x|\ge m_k}\int_{\R^3} e^{\lambda_{f_0}} f_{n_{m_k}}\,dx\,dv
   +\bigg|\int_{|x|\le m_k}\int_{\R^3} e^{\lambda_{f_0}} f_{n_{m_k}}\,dx\,dv
   -\theta_0\bigg|
   \\ & \le & \int_{|x|\ge m_k}\int_{\R^3} e^{\lambda_{f_0}} f_{n_{m_k}}\,dx\,dv
   +\eps_{m_k}\to 0,\quad k\to\infty,  
\end{eqnarray*} 
which yields the assertion (\ref{ropla}) along the specified subsequence.  
Case 2: There is a $\delta_0>0$ such that 
\[ \int_{|x|\ge m}\int_{\R^3} e^{\lambda_{f_0}} f_{n_m}\,dx\,dv\ge\delta_0,
   \quad m\ge m_0. \] 
Then (\ref{cbetabd}) shows that  
\[ \int_{|x|\ge m}\int_{\R^3}\sqrt{1+|v|^2} f_{n_m}\,dx\,dv\ge\delta_1,
   \quad m\ge m_0, \]  
where $\delta_1=c_\beta^{-1}\delta_0>0$. Since every term on the left-hand side 
is bounded from above by the mass, there is a further subsequence (not relabelled) 
such that 
\[ \eps_0=\lim_{m\to\infty}\int_{|x|\ge m}\int_{\R^3}\sqrt{1+|v|^2} f_{n_m}\,dx\,dv
   \in [\delta_1, M] \] 
does exist. Now $\eps_0=M$ is impossible, since then 
\begin{eqnarray*} 
   \int_{|x|\le R_0}\int_{\R^3}\sqrt{1+|v|^2} f_{n_m}\,dx\,dv
   & \le & \int_{|x|\le m}\int_{\R^3}\sqrt{1+|v|^2} f_{n_m}\,dx\,dv
   \\ & = & M-\int_{|x|\ge m}\int_{\R^3}\sqrt{1+|v|^2} f_{n_m}\,dx\,dv\to M-\eps_0=0
\end{eqnarray*}  
as $m\to\infty$, but 
\[ \lim_{m\to\infty}\int_{|x|\le R_0}\int_{\R^3}\sqrt{1+|v|^2}\,f_{n_m}\,dx\,dv
   =\int_{|x|\le R_0}\int_{\R^3}\sqrt{1+|v|^2}\,f_0\,dx\,dv=M \] 
due to (\ref{eram}) and the support properties of $f_0$. 
Altogether, we have seen that $\eps_0\in ]0, M[$. 
Since also ${(f_{n_m})}_{m\in\N}$ is a minimizing sequence and $R_m=m$ 
will satisfy (\ref{Rncond2}) for large $m$, Corollary \ref{SBECreal} yields a contradiction. 
{\hfill$\Box$}\bigskip 

The following theorem summarizes the foregoing results. 
Up to the compact support property of $f_0$ in $v$ that is shown later 
(up to now we only know (\ref{tiha})), it contains all the statements 
of the main theorem, Theorem \ref{mainthm}. 

\begin{theorem}\label{f0exists} 
Let $k\in ]0, 2]$. Then (CBEC) is verified.  
Suppose that $f_0\in L^{1+1/k}(\R^3\times\R^3)$ 
is the limiting function as obtained in Lemma \ref{f0props}.
Then $f_0$ is a minimizer, i.e., $f_0\in {\cal A}$ and ${\cal D}(f_0)=I$. 
Furthermore, $f_0\le 1$, 
\begin{equation}\label{tiha} 
   \int_{\R^3}\int_{|v|\ge P+1}\sqrt{1+|v|^2}\,f_0\,dx\,dv
   \le\frac{2}{P^{1/4}},\quad P\ge P_0+1,
\end{equation} 
and there is $R_0>0$ such that $f_0(x, v)=0$ for a.e.~$|x|\ge R_0$ and $v\in\R^3$. 
\end{theorem}
{\bf Proof\,:} From Remark \ref{sbech}(b) we know that (CBEC) holds, 
and by Lemma \ref{f0incalA} we have $f_0\in {\cal A}$. 
We write the functional as 
\[ {\cal D}(f)=\Psi(f)-\int_{\R^3}\int_{\R^3} e^{\lambda_f} f\,dx\,dv,
   \quad\Psi(f)=\int_{\R^3}\int_{\R^3} e^{\lambda_f} \chi(f)\,dx\,dv, \] 
and consider ${\cal D}(f_n)$, disregarding the fact that we would have 
to pass to subsequences in what follows. To deal with the term $\Psi(f_n)$, 
we are going to use Lemma \ref{phikonv}, which asserts that $\Psi$ is convex.
Since $(f_n)\subset {\cal A}$ and also ${\cal A}$ is convex,
it follows from Mazur's Theorem \cite[Thm.~3.9]{dacor} and (\ref{weakc}) that
$\Psi(f_0)\le\liminf_{n\to\infty}\Psi(f_n)$. Hence by Lemma \ref{2ndpartconv} 
\begin{eqnarray*} 
   {\cal D}(f_0) & = & \Psi(f_0)-\int_{\R^3}\int_{\R^3} e^{\lambda_0} f_0\,dx\,dv
   \\ & \le & \liminf_{n\to\infty}\bigg(\Psi(f_n)-\int_{\R^3}\int_{\R^3} e^{\lambda_{f_n}} f_n\,dx\,dv\bigg)
   =\liminf_{n\to\infty} {\cal D}(f_n)=I. 
\end{eqnarray*} 
Since $f_0\in {\cal A}$, also $I\le {\cal D}(f_0)$, which yields ${\cal D}(f_0)=I$. 
For the stated additional properties of $f_0$, see Lemma \ref{f0props} and Lemma \ref{f0cpt}. 
{\hfill$\Box$}\bigskip

%%%%%%%%%%%%%%%%%%%%%%%%%%%%%%%%%%%%%%%%%%%%%%%%%%%%%%%%%%%%%%%%%%%%%%%%%%%%%%%%%%%%%%%%%

\section{A variational inequality and compact support in $v$-space}
\label{furthprop_sect} 
\setcounter{equation}{0}

From Lemma \ref{f0cpt} we already know that $f_0(x, v)=0$ for a.e.~$|x|\ge R_0$ and $v\in\R^3$. 
Also $0\le f_0\le 1$ a.e.~by construction. 
Denote $\rho_0=\rho_{f_0}$, $m_0=m_{f_0}$, $\lambda_0=\lambda_{f_0}$, etc. 

\begin{lemma}\label{formlin} If $g\in {\cal A}$, then  
\begin{equation}\label{DHdiffb} 
   \int\int\Big[
   e^{\lambda_0} (\chi'(f_0)-1) (g-f_0)
   +e^{3\lambda_0}\,(\chi(f_0)-f_0)\,\frac{m_g-m_0}{r}\Big]\,dx\,dv\ge 0
\end{equation} 
\end{lemma} 
{\bf Proof\,:} Let $g\in {\cal A}$. From the convexity of ${\cal A}$ 
and the minimality of $f_0$ we can infer that $\varphi(t)={\cal D}(tg+(1-t)f_0)
\ge {\cal D}(f_0)=\varphi(0)$ for $t\in [0, 1]$. Therefore by definition of ${\cal D}$ 
\begin{equation}\label{hander} 
   \int_{\R^3}\int_{\R^3} e^{\lambda_0}\,(\chi(f_0)-f_0)\,dx\,dv
   \le\int_{\R^3}\int_{\R^3} e^{\lambda_{tg+(1-t)f_0}}
   \,\Big(\chi(tg+(1-t)f_0)-f_0-t(g-f_0)\Big)\,dx\,dv
\end{equation}   
for $t\in [0, 1]$. To begin with, we recall from (\ref{cbetabd}) that 
$1\le e^{\lambda_g}\le c_\beta$ for all $g\in {\cal A}$. Also 
\[ \Big|\frac{m_g(r)}{r}-\frac{m_0(r)}{r}\Big|
   =\frac{4\pi}{r}\,\Big|\int_0^r s^2\,(\rho_g(s)-\rho_0(s))\,ds\Big|. \] 
Hence if $0<r\le 1$, then 
\[ \Big|\frac{m_g(r)}{r}-\frac{m_0(r)}{r}\Big|
   \le\frac{4\pi}{r}\,\Big(\int_0^r s^2\,(\rho_g(s)+\rho_0(s))\,ds\Big)
   \le\frac{8\pi\sigma_0}{3r}\,r^3\le\frac{8\pi}{3}\,\sigma_0. \] 
On the other hand, if $r\ge 1$, then 
\[ \Big|\frac{m_g(r)}{r}-\frac{m_0(r)}{r}\Big|
   \le\frac{4\pi}{r}\,\int_0^r s^2\,(\rho_g(s)+\rho_0(s))\,ds
   \le 2M. \] 
Thus we get 
\begin{equation}\label{haste} 
   2e^{2\lambda_0(r)}\Big|\frac{m_g(r)-m_0(r)}{r}\Big|
   \le 2c_g^2\,\Big(\frac{8\pi}{3}\,\sigma_0+2M\Big)=:C_0
\end{equation}  
in any case. Henceforth we restrict the range of $t$ 
to $t\in [0, 1]\cap [0, \frac{1}{2C_0}]$. Since $m_f$ is linear in $f$, 
we can write pointwise a.e. 
\begin{eqnarray*} 
   e^{-2\lambda_{tg+(1-t)f_0}} & = & 1-\frac{2m_{tg+(1-t)f_0}}{r}
   =1-\frac{2m_0}{r}-2t\,\frac{m_g-m_0}{r}
   =e^{-2\lambda_0}-2t\,\frac{m_g-m_0}{r}
   \\ & = & e^{-2\lambda_0}\Big(1-2t\,e^{2\lambda_0}\frac{m_g-m_0}{r}\Big). 
\end{eqnarray*}  
In view of (\ref{haste}) and $0\le t\le\frac{1}{2C_0}$, it follows that 
\[ e^{\lambda_{tg+(1-t)f_0}} 
   =e^{\lambda_0}\frac{1}{(1-2t\,e^{2\lambda_0}\frac{m_g-m_0}{r})^{1/2}}. \] 
We have 
\[ \frac{1}{(1-x)^{1/2}}=1+\frac{1}{2}\,x+\psi(x)x^2,
   \quad |\psi(x)|\le C,\quad x\in\Big[0, \frac{1}{2}\Big]. \] 
Therefore we may continue     
\begin{eqnarray*} 
   e^{\lambda_{tg+(1-t)f_0}} 
   & = & e^{\lambda_0}\Big(1+t\,e^{2\lambda_0}\frac{m_g-m_0}{r}
   +\psi\Big(2t\,e^{2\lambda_0}\frac{m_g-m_0}{r}\Big)
   \Big[2t\,e^{2\lambda_0}\frac{m_g-m_0}{r}\Big]^2\Big)
   \\ & = & e^{\lambda_0}\Big(1+t\,e^{2\lambda_0}\frac{m_g-m_0}{r}+t^2\Psi\Big), 
\end{eqnarray*} 
where $|\Psi(t, r)|\le C$ is uniformly bounded. Next, 
\begin{eqnarray*} 
   \chi(tg+(1-t)f_0) & = & \chi(f_0)+t\chi'(f_0)(g-f_0)
   +\int_0^{t(g-f_0)}\Big(\chi'(f_0+s)-\chi'(f_0)\Big)\,ds
   \\ & = & \chi(f_0)+t\chi'(f_0)(g-f_0)+t(g-f_0)\,\Phi,  
   \\ \Phi(t, x, v) & = & \frac{1}{t}\int_0^t
   \Big(\chi'(f_0(x, v)+\tau(g(x, v)-f_0(x, v)))-\chi'(f_0(x, v))\Big)\,d\tau.  
\end{eqnarray*} 
Going back to (\ref{hander}) and using the previous expansions, we get 
for $t\in ]0, \frac{1}{2C_0}]$, after some simplification, 
\begin{eqnarray*} 
   0 & \le & \int_{\R^3}\int_{\R^3}
   e^{\lambda_0}\bigg(
   \Big[(\chi'(f_0)-1)(g-f_0)+e^{2\lambda_0}(\chi(f_0)-f_0)\,\frac{m_g-m_0}{r}
        +(g-f_0)\Phi\Big] 
        \\ & & \hspace{6em} +\,t\Big[e^{2\lambda_0}\frac{m_g-m_0}{r}\,(g-f_0)
        (\chi'(f_0)-1+\Phi)+(\chi(f_0)-f_0)\Psi\Big] 
        \\ & & \hspace{6em} +\,t^2\Big[\chi'(f_0)(g-f_0)\Psi+(g-f_0)\Psi\,\Phi
        -(g-f_0)\Psi\Big]\bigg)\,dx\,dv.    
\end{eqnarray*} 
It remains to show that the error terms tend to zero as $t\to 0$. 
For, we will only prove that 
\begin{equation}\label{Ito0} 
   I(t)=\int_{\R^3}\int_{\R^3} e^{\lambda_0}(g-f_0)\Phi\,dx\,dv\to 0,
   \quad t\to 0,
\end{equation}   
since the other terms are similar or easier. For a.e.~$(x, v)\in\R^3\times\R^3$ 
we have $\Phi(t, x, v)\to 0$ as $t\to 0$, since $\chi'$ is continuous on $[0, \infty[$, 
$0\le f_0(x, v)\le 1$ a.e.~and $0\le g(x, v)<\infty$ a.e.; 
also note that $f_0(x, v)+\tau(g(x, v)-f_0(x, v))=(1-\tau)f_0(x, v)+\tau g(x, v)\ge 0$.  
No matter if $k\le 1$ or $k\ge 1$, there is a constant $C_k>0$ 
such that $\chi'(s_1+s_2)=(s_1+s_2)^{1/k}\le C_k(s_1^{1/k}+s_2^{1/k})$ 
for $s_1, s_2\ge 0$. Therefore pointwise a.e. 
\begin{eqnarray*} 
   |\Phi(t, x, v)| & \le & \frac{1}{t}\int_0^t 
   \Big(C_k(1-\tau)^{1/k} f_0(x, v)^{1/k}+C_k\tau^{1/k} g(x, v)^{1/k}
   +f_0(x, v)^{1/k}\Big)\,d\tau
   \\ & \le & (C_k+1)(f_0(x, v)^{1/k}+g(x, v)^{1/k}),      
\end{eqnarray*}   
and hence  
\[  e^{\lambda_0}\,|g(x, v)-f_0(x, v)|\,|\Phi(t, x, v)|
    \le c_\beta (C_k+1)(f_0(x, v)+g(x, v))\,(f_0(x, v)^{1/k}+g(x, v)^{1/k}). \] 
Since $f_0, g\in L^{1+1/k}(\R^3\times\R^3)$, the right-hand side is integrable, 
by H\"older's inequality. Therefore the dominated convergence theorem 
shows that (\ref{Ito0}) is valid. 
{\hfill$\Box$}\bigskip

\begin{remark}\label{tipe}{\rm If $\rho_g=\rho_0$ and $g\geq 0$, then $g\in {\cal A}$ is for free, 
and due to $m_g=m_0$ the second term in (\ref{DHdiffb}) drops out. Thus we have 
\begin{equation}\label{DHdiff2}  
   \rho_g=\rho_0\quad\Longrightarrow\quad
   \int\int_{{\rm supp}\,f_0} e^{\lambda_0} (1-\chi'(f_0)) (f_0-g)\,dx\,dv
   \ge\int\int_{({\rm supp}\,f_0)^c} e^{\lambda_0} g\,dx\,dv. 
\end{equation} 
}
\end{remark} 

This may be used to prove that $f_0$ has compact support in $v$. 
Taking together Theorem \ref{f0exists} and Corollary \ref{vsuppcpt}, 
this in particular completes the proof of Theorem \ref{mainthm}. 
\begin{cor}\label{vsuppcpt} 
In the setting of Theorem \ref{f0exists} the following holds for the minimizer $f_0$. 
There is a $S_0>0$ such that $f_0(x, v)=0$ for a.e.~$x\in\R^3$ and $|v|\ge S_0$. 
\end{cor}
{\bf Proof\,:} We consider the boxes 
\begin{equation}\label{BhatB} 
   B=\{(x, v): |x|\in [r_1, r_2], |v|\in [s_1, s_2]\},
   \quad\hat{B}=\{(x, v): |x|\in [r_1, r_2], |v|\in [\hat{s}_1, \hat{s}_2]\},
\end{equation} 
for $r_2>r_1>0$, $s_2>s_1>0$ and $\hat{s}_2>\hat{s}_1>0$ such that $B\cap\hat{B}=\emptyset$ 
and $B$ is below $\hat{B}$. Define 
\[ g(x, v)=\left\{\begin{array}{c@{\quad:\quad}c}
   f_0(x, v)+b(x, v) & (x, v)\in B 
   \\[1ex]
   \frac{1}{2}f_0(x, v) & (x, v)\in\hat{B} 
   \\[1ex] 
   f_0(x, v) & (x, v)\in\R^6\setminus (B\cup\hat{B}) 
   \end{array}\right. , \] 
where 
\[ b(x, v)=\bigg(\frac{1}{2}\int_{\hat{s}_1\le |\tilde{v}|\le\hat{s}_2}
   \sqrt{1+|\tilde{v}|^2}\,f_0(x, \tilde{v})\,d\tilde{v}\bigg)\,
   \frac{1}{4\pi(s_2-s_1)}\,\frac{1}{|v|^2\sqrt{1+|v|^2}}. \] 
Then $\rho_g(r)=\rho_0(r)$ for $r\in [0, \infty[\setminus [r_1, r_2]$. 
If $r\in [r_1, r_2]$ and since $B$ is below $\hat{B}$, we get 
\begin{eqnarray*}
   \rho_g(r) & = & \int_{0\le |v|\le s_1}\sqrt{1+|v|^2}\,f_0(x, v)\,dv
   +\int_{s_1\le |v|\le s_2}\sqrt{1+|v|^2}\,(f_0(x, v)+b(x, v))\,dv
   \\ & & +\int_{s_2\le |v|\le\hat{s}_1}\sqrt{1+|v|^2}\,f_0(x, v)\,dv
   +\frac{1}{2}\int_{\hat{s}_1\le |v|\le\hat{s}_2}\sqrt{1+|v|^2}\,f_0(x, v)\,dv
   \\ & & +\int_{|v|\ge\hat{s}_2}\sqrt{1+|v|^2}\,f_0(x, v)\,dv
   \\ & = & \rho_{f_0}(r)+\int_{s_1\le |v|\le s_2}\sqrt{1+|v|^2}\,b(x, v)\,dv
   -\frac{1}{2}\int_{\hat{s}_1\le |v|\le\hat{s}_2}\sqrt{1+|v|^2}\,f_0(x, v)\,dv
   \\ & = & \rho_{f_0}(r)
   +\bigg(\frac{1}{2}\int_{\hat{s}_1\le |\tilde{v}|\le\hat{s}_2}
   \sqrt{1+|\tilde{v}|^2}\,f_0(x, \tilde{v})\,d\tilde{v}\bigg)
   \,\frac{4\pi}{4\pi(s_2-s_1)}\int_{s_1}^{s_2}\frac{u^2\sqrt{1+u^2}}{u^2\sqrt{1+u^2}}\,\,du
   \\ & & -\,\frac{1}{2}\int_{\hat{s}_1\le |v|\le\hat{s}_2}\sqrt{1+|v|^2}\,f_0(x, v)\,dv
   \\ & = & \rho_0(r). 
\end{eqnarray*}  
Thus we have seen that $\rho_g=\rho_0$ and (\ref{DHdiff2}) applies. It follows that  
\[ \int\int_{{\rm supp}\,f_0} e^{\lambda_0} (1-\chi'(f_0)) (f_0-g)\,dx\,dv
   \ge\int\int_{({\rm supp}\,f_0)^c} e^{\lambda_0} g\,dx\,dv, \] 
which leads to 
\begin{eqnarray*} 
   \lefteqn{\frac{1}{2}\int\int_{\hat{B}} e^{\lambda_0} (1-\chi'(f_0)) f_0\,dx\,dv} 
   \\ & = & \frac{1}{2}\int\int_{{\rm supp}\,f_0\cap\hat{B}} e^{\lambda_0} (1-\chi'(f_0)) f_0\,dx\,dv
   \\ & \ge & \int\int_{{\rm supp}\,f_0\cap B} e^{\lambda_0} (1-\chi'(f_0))\,b\,dx\,dv
   +\int\int_{({\rm supp}\,f_0)^c\cap B} e^{\lambda_0}\,b\,dx\,dv
   \\ & = & \int\int_B e^{\lambda_0} (1-\chi'(f_0))\,b\,dx\,dv
   \\ & = & \int_{r_1\le |x|\le r_2} dx\,e^{\lambda_0(r)}\int_{s_1\le |v|\le s_2} dv
   \,(1-\chi'(f_0(x, v)))\,\bigg(\frac{1}{2}\int_{\hat{s}_1\le |\tilde{v}|\le\hat{s}_2}
   \sqrt{1+|\tilde{v}|^2}\,f_0(x, \tilde{v})\,d\tilde{v}\bigg)\,
   \\ & & \hspace{12em}\times\frac{1}{4\pi(s_2-s_1)}\,\frac{1}{|v|^2\sqrt{1+|v|^2}}
   \\ & = & \frac{1}{2}\int_{r_1\le |x|\le r_2} dx\,e^{\lambda_0(r)}
   \int_{\hat{s}_1\le |v|\le\hat{s}_2} dv\,\sqrt{1+|v|^2}\,f_0(x, v)
   \int_{s_1\le |\tilde{v}|\le s_2} d\tilde{v}
   \,(1-\chi'(f_0(x, \tilde{v})))\,
   \\ & & \hspace{12em}\times\frac{1}{4\pi(s_2-s_1)}\,\frac{1}{|\tilde{v}|^2\sqrt{1+|\tilde{v}|^2}}. 
\end{eqnarray*} 
This can be re-expressed as 
\begin{eqnarray}\label{tggp}  
   \lefteqn{\int\int_{\hat{B}} dx\,dv\,e^{\lambda_0(r)} (1-\chi'(f_0(x, v))) f_0(x, v)} 
   \nonumber
   \\ & \ge & \frac{1}{4\pi(s_2-s_1)}\int\int_{\hat{B}}\,dx\,dv\,e^{\lambda_0(r)}
   \,\sqrt{1+|v|^2}\,f_0(x, v)\\
& & \hspace{12em}  \int_{s_1\le |\tilde{v}|\le s_2} d\tilde{v}
   \,(1-\chi'(f_0(x, \tilde{v})))\,\frac{1}{|\tilde{v}|^2\sqrt{1+|\tilde{v}|^2}}.  
   \nonumber \\ & & 
\end{eqnarray} 
This relation is verified for all boxes $\hat{B}$ above $B$ as in (\ref{BhatB}) 
that do not intersect. Next recall from $f_0\in {\cal A}$ and Theorem \ref{f0exists} 
that $\rho_0(r)\le\sigma_0$ and $f_0\le 1$ a.e. Take 
\[ N\ge\max\Big\{\sqrt[3]{\frac{3\sigma_0}{2\pi}}, \sqrt[3]{\frac{3\sigma_0}{2\pi}\,2^k}\Big\}. \] 
To begin with, for a.e.~$x\in\R^3$, 
\begin{eqnarray}\label{expim}  
   \lefteqn{\int_{|\tilde{v}|\le N} d\tilde{v}
   \,(1-\chi'(f_0(x, \tilde{v})))\,\frac{1}{|\tilde{v}|^2\sqrt{1+|\tilde{v}|^2}}}
   \nonumber 
   \\ & \ge & \frac{1}{N^2\sqrt{1+N^2}}\int_{|\tilde{v}|\le N} d\tilde{v}
   \,(1-f_0(x, \tilde{v})^{1/k})
   \nonumber
   \\ & = & \frac{1}{N^2\sqrt{1+N^2}}\bigg(\frac{4\pi}{3}\,N^3
   -\int_{|\tilde{v}|\le N} d\tilde{v}\,f_0(x, \tilde{v})^{1/k}\bigg). 
\end{eqnarray} 
If $k\le 1$, then $f_0^{1/k}\le f_0$, and hence 
\begin{eqnarray*}
   \int_{|\tilde{v}|\le N} d\tilde{v}
   \,(1-\chi'(f_0(x, \tilde{v})))\,\frac{1}{|\tilde{v}|^2\sqrt{1+|\tilde{v}|^2}}
   & \ge & \frac{1}{N^2\sqrt{1+N^2}}\bigg(\frac{4\pi}{3}\,N^3-\rho_{f_0}(r)\bigg)
   \\ & \ge & \frac{1}{N^2\sqrt{1+N^2}}\bigg(\frac{4\pi}{3}\,N^3-\sigma_0\bigg)
   \\ & \ge & \frac{2\pi}{3}\,\frac{N}{\sqrt{1+N^2}}. 
\end{eqnarray*}
On the other hand, if $k\ge 1$ and $f_0(x, \tilde{v})^{1/k}\ge 1/2$, then 
\[ 2^{k-1} f_0(x, \tilde{v})f_0(x, \tilde{v})^{-1/k}
   =(2 f_0(x, \tilde{v})^{1/k})^{k-1}\ge 1 \] 
yields $f_0(x, \tilde{v})^{1/k}\le 2^{k-1} f_0(x, \tilde{v})$. It follows that 
\begin{eqnarray*} 
   \int_{|\tilde{v}|\le N} d\tilde{v}\,f_0(x, \tilde{v})^{1/k}
   & = & \int_{\{|\tilde{v}|\le N,\,f_0(x, \tilde{v})^{1/k}\le 1/2\}} d\tilde{v}\,f_0(x, \tilde{v})^{1/k}
   +\int_{\{|\tilde{v}|\le N,\,f_0(x, \tilde{v})^{1/k}\ge 1/2\}} d\tilde{v}\,f_0(x, \tilde{v})^{1/k}
   \\ & \le & \frac{1}{2}\,\frac{4\pi}{3}\,N^3
   +2^{k-1}\int_{\{|\tilde{v}|\le N,\,f_0(x, \tilde{v})^{1/k}\ge 1/2\}} d\tilde{v}\,f_0(x, \tilde{v})
   \\ & \le & \frac{2\pi}{3}\,N^3+2^{k-1}\rho_0(r). 
\end{eqnarray*} 
Thus in the case $k\ge 1$ we deduce from (\ref{expim}) and the choice of $N$ that 
\begin{eqnarray*}  
   \lefteqn{\int_{|\tilde{v}|\le N} d\tilde{v}
   \,(1-\chi'(f_0(x, \tilde{v})))\,\frac{1}{|\tilde{v}|^2\sqrt{1+|\tilde{v}|^2}}}
   \\ & \ge & \frac{1}{N^2\sqrt{1+N^2}}\bigg(\frac{2\pi}{3}\,N^3-2^{k-1}\rho_{f_0}(r)\bigg) 
   \\ & \ge & \frac{1}{N^2\sqrt{1+N^2}}\bigg(\frac{2\pi}{3}\,N^3-2^{k-1}\sigma_0\bigg)
   \\ & \ge & \frac{\pi}{3}\,\frac{N}{\sqrt{1+N^2}}.  
\end{eqnarray*} 
Therefore we have seen that in all cases the estimate 
\[ \int_{|\tilde{v}|\le N} d\tilde{v}
   \,(1-\chi'(f_0(x, \tilde{v})))\,\frac{1}{|\tilde{v}|^2\sqrt{1+|\tilde{v}|^2}}
   \ge\frac{\pi}{3}\,\frac{N}{\sqrt{1+N^2}} \] 
is verified. Thus taking $s_1=0$ and $s_2=N$ in (\ref{tggp}), 
we therefore obtain for any box $\hat{B}$ that lies above the line $\{|v|=N\}$ the estimate
\begin{eqnarray*} 
   \lefteqn{12\sqrt{1+N^2}\int\int_{\hat{B}} dx\,dv\,e^{\lambda_0(r)} (1-\chi'(f_0(x, v))) f_0(x, v)} 
   \\ & \ge & \int\int_{\hat{B}}\,dx\,dv\,e^{\lambda_0(r)}
   \,\sqrt{1+|v|^2}\,f_0(x, v).  
\end{eqnarray*} 
Since this holds for all such boxes $\hat{B}$, we must have pointwise a.e.~in the set 
$\{(x, v)\in\R^3\times\R^3: |v|\ge N\}$ the relation 
\[ \sqrt{1+|v|^2}\,f_0(x, v)\le 12\sqrt{1+N^2}\,(1-\chi'(f_0(x, v))) f_0(x, v)
   \le 12\sqrt{1+N^2}\,f_0(x, v), \] 
and we can in turn infer that $|v|\le 12\sqrt{1+N^2}$ on ${\rm supp}\,f_0$. 
{\hfill$\Box$}\bigskip

%%%%%%%%%%%%%%%%%%%%%%%%%%%%%%%%%%%%%%%%%%%%%%%%%%%%%%%%%%%%%%%%%%%%%%%%%%%%%%%    

\section{When is the minimizer a static solution?}
\label{isminstatic_sect}
\setcounter{equation}{0}

We recall some results from \cite{AK1} (with $\alpha=1$ there), 
which are stated in the notation that we are using in the present paper. 
Let $r>0$ be fixed and put
\[ \Xi_r=\Big\{\psi=\psi(w, \ell)\in L^{1+1/k}(\R\times [0, \infty[): 
   \psi\ge 0, \psi(w, \ell)=0\,\,\mbox{for}\,\,\sqrt{1+w^2+\ell^2/r^2}\ge\sqrt{1+S_0^2}\Big\}. \] 
As compared to \cite{AK1}, the number `$2$' in the definition of $\Xi_r$ 
has been replaced by $\sqrt{1+S_0^2}$, and as remarked in \cite{AK1}, 
the results still hold with such a modification. Now we consider the functionals 
\begin{equation}\label{HF-def} 
   \hat{H}(\psi)=\int_0^\infty d\ell^2\int_{\R} dw\,(\chi(\psi)-\psi),
   \quad\hat{F}(\psi, r, a)=\int_0^\infty d\ell^2\int_{\R} dw\,
   \sqrt{1+w^2+\ell^2/r^2}\,\psi-\frac{r^2}{\pi}\,a,
\end{equation} 
for $r>0$, functions $\psi=\psi(w, \ell)\in\Xi_r$ and $a>0$; 
here $d\ell^2=2\ell\,d\ell$.  
The next lemma is \cite[Lemma 5.1]{AK1}. 

\begin{lemma}\label{wibirk} For $r>0$ and $a>0$ 
there is a unique function $\psi=\psi(w, \ell; r, a)\in\Xi_r$ such that 
\begin{equation}\label{lslem} 
   \inf_{\psi\in\Xi_r}\,\{\hat{H}(\psi): \hat{F}(\psi, r, a)=0\}
   =\hat{H}(\psi(\cdot, \cdot; r, a))
\end{equation}  
and 
\begin{equation}\label{F0} 
   \hat{F}(\psi(\cdot, \cdot; r, a), r, a)=0.
\end{equation}
For $r>0$ and $a=0$ this unique function is $\psi(w, \ell; r, 0)=0$. 
\end{lemma} 

\begin{cor}\label{AK1cor} We have 
\begin{equation}\label{trythi} 
   \psi(\cdot, \cdot; r, \rho_0(r))=f_0(r, \cdot, \cdot)
\end{equation}  
for a.e.~$r\ge 0$. 
\end{cor}
{\bf Proof\,:} Since 
\[ 4\pi^2\int_0^\infty dr\int_{\R} dw\int_0^\infty d\ell^2\, f_0^{1+1/k}
   =\int_{\R^3}\int_{\R^3} f_0^{1+1/k}\,dx\,dv<\infty, \] 
we have $\int_{\R} dw\int_0^\infty d\ell^2\, f_0(r, \cdot, \cdot)^{1+1/k}<\infty$ 
for a.e.~$r\ge 0$, so that $0\le f_0(r, \cdot, \cdot)\in L^{1+1/k}(\R\times [0, \infty[)$. 
Also, by Corollary \ref{vsuppcpt}, the support condition in $\Xi_r$ is satisfied, 
and hence $f_0(r, \cdot, \cdot)\in\Xi_r$. Next observe that 
\[ \hat{F}(f_0(r, \cdot, \cdot), r, \rho_0(r))
   =\int_0^\infty d\ell^2\int_{\R} dw\,
   \sqrt{1+w^2+\ell^2/r^2}\,f_0(r, w, \ell)-\frac{r^2}{\pi}\,\rho_0(r)=0 \] 
due to $dv=\frac{\pi}{r^2}\,d\ell^2\,dw$. It follows from the 
minimizing property of the function $\psi$ that 	
\begin{equation}\label{sas} 
   \hat{H}(\psi(\cdot, \cdot; r, \rho_0(r)))\le\hat{H}(f_0(r, \cdot, \cdot)).
\end{equation}   
Let the spherically symmetric function $g\ge 0$ be defined by 
\[ g(r, w, \ell)=\psi(w, \ell; r, \rho_0(r))\ge 0. \] 
Then $\rho_g=\rho_0\le\sigma_0$ by the constraint $\hat{F}=0$. 
In particular, $g$ has total mass $M$. Next, (\ref{sas}) says that 
\[ \int_0^\infty d\ell^2\int_{\R} dw\,(\chi(g(r, \cdot, \cdot))-g(r, \cdot, \cdot))
   \le\int_0^\infty d\ell^2\int_{\R} dw\,(\chi(f_0(r, \cdot, \cdot))-f_0(r, \cdot, \cdot)), \] 
and upon integration $\int_0^\infty dr$ it follows that 
\begin{eqnarray*} 
   \frac{k}{k+1}\int_{\R^3}\int_{\R^3} g^{1+1/k}\,dx\,dv
   & \le & \int_{\R^3}\int_{\R^3} g\,dx\,dv
   +\frac{k}{k+1}\int_{\R^3}\int_{\R^3} f_0^{1+1/k}\,dx\,dv
   \\ & \le & M+\frac{k}{k+1}\|f_0\|_{L^{1+1/k}(\R^3\times\R^3)}^{1+1/k}<\infty.
\end{eqnarray*} 
As a consequence, $g\in L^{1+1/k}(\R^3\times\R^3)$, which shows that $g\in {\cal A}$. 
Since $f_0$ is a minimizer, and $m_g=m_0$ yields $\lambda_g=\lambda_0$,   
\begin{eqnarray*} 
   0 & \ge & {\cal D}(f_0)-{\cal D}(g)
   =4\pi^2\int_0^\infty dr\,e^{\lambda_0}
   \int_{\R} dw\int_0^\infty d\ell^2\,[(\chi(f_0)-f_0)-(\chi(g)-g)]
   \\ & = & 4\pi^2\int_0^\infty dr\,e^{\lambda_0} [\hat{H}(f_0(r, \cdot, \cdot))
   -\hat{H}(g(r, \cdot, \cdot))]\ge 0, 
\end{eqnarray*} 
where in the last step (\ref{sas})  has been used. Hence 
$\hat{H}(f_0(r, \cdot, \cdot))=\hat{H}(g(r, \cdot, \cdot))$ 
for a.e.~$r\ge 0$, and from the uniqueness assertion in Lemma \ref{wibirk} 
we infer that (\ref{trythi}) holds. 
{\hfill$\Box$}\bigskip

Now we can follow the steps in \cite{AK1} and arrive at \cite[equ.~(5.15)]{AK1} 
to deduce that 
\begin{equation}\label{fastform} 
   f_0(r, w, \ell)
   =\Big(1-\tilde{\eps}(r)\sqrt{1+w^2+\ell^2/r^2}\Big)_+^k 
\end{equation}
for $\tilde{\eps}(r)=(G')^{-1}(-4\pi\kappa\rho_0(r))$, 
where $\kappa=\frac{k+1}{16\pi^2}$ and 
\[ G(\eps)=\int_0^\infty\xi^2
   \,{(1-\eps\sqrt{1+\xi^2})}_+^{k+1}\,d\xi,\quad\eps\in\R; \] 
see \cite[Section 7.1]{AK1} for the properties of $G$ and its derivatives.    
It is important to note that the fact that $\rho_\ast$ from \cite{AK1} satisfies 
the (reformulated) Euler-Lagrange equation enters only in the end, 
where it is shown that the function $\tilde{\eps}$ 
is proportional to $e^{\mu_\ast}$. 
\smallskip

This allows the following reformulation of (\ref{DHdiffb}). 

\begin{cor}\label{reform} If $g\in {\cal A}$, then 
\begin{equation}\label{firel}   
   \int_0^\infty dr\int_0^\infty d\ell^2\int_{\R} dw
   \,e^{\lambda_0}\bigg[(\chi'(f_0)-1)(g-f_0)
   +e^{2\lambda_0}\,f_0\frac{d}{dw}\,[(\chi'(f_0)-1)w] 
   \,\frac{m_g-m_0}{r}\bigg]\ge 0.   
\end{equation} 
\end{cor} 
{\bf Proof\,:} From Lemma \ref{f0cpt} and Corollary \ref{vsuppcpt} 
we know that $f_0(r, w, \ell)=0$ for $r\ge R_0$ or if $1+w^2+\ell^2/r^2\ge 1+S_0^2$. 
If $(r, \ell)$ are fixed, then (\ref{fastform}) implies that the support of $f_0(r, \cdot, \ell)$ 
is given by the (symmetric) interval $I_{r, \ell}=\{w\in\R: |w|^2\le 1/\tilde{\eps}(r)^2-1-\ell^2/r^2\}$, 
and $f_0(r, w, \ell)=(1-\tilde{\eps}(r)\sqrt{1+w^2+\ell^2/r^2})^k$ on $I_{r, \ell}$, 
which is a differentiable function in $w$. Therefore we may integrate by parts in $w$ 
to obtain 
\begin{eqnarray*}
   \int_{\R^3} (\chi (f_0)-f_0)\,dv 
   & = & \frac{\pi}{r^2}\int_0^\infty d\ell^2\int_{\R} dw\,(\chi(f_0)-f_0) 
   =\frac{\pi}{r^2}\int_0^\infty d\ell^2\int_{I_{r, \ell}} dw\,(\chi(f_0)-f_0)
   \,\frac{d}{dw}\,w 
   \\ & = & -\frac{\pi}{r^2}\int_0^\infty d\ell^2\int_{I_{r, \ell}} dw\,
   (\chi'(f_0)-1)\,\frac{\partial f_0}{\partial w}\,w
   \\ & = & \frac{\pi}{r^2}\int_0^\infty d\ell^2\int_{I_{r, \ell}} dw\,
   \frac{d}{dw}\,[(\chi'(f_0)-1)w]\,f_0
   \\ & = & \frac{\pi}{r^2}\int_0^\infty d\beta\int_{\R} dw\,
   \frac{d}{dw}\,[(\chi'(f_0)-1)w]\,f_0. 
\end{eqnarray*}
Hence (\ref{firel}) follows from (\ref{DHdiffb}) 
and $dx\,dv=4\pi^2 dr\,d\ell^2\,dw$. 
{\hfill$\Box$}\bigskip 

The relation (\ref{firel}) can also be restated in a slimmer form. 

\begin{cor} Let $I_{r, \ell}=\{w\in\R: |w|^2\le 1/\tilde{\eps}(r)^2-1-\ell^2/r^2\}$ and 
\[ U(r)=-\left\{\begin{array}{c@{\quad:\quad}l} e^{-\mu_0(r)}\tilde{\eps}(r) & r\in [0, R_0] 
   \\ 0 & r\in ]R_0, \infty[\end{array}\right. . \] 
If $g\in {\cal A}$, then 
\[ \int_0^\infty\frac{d}{dr}\Big(e^{\lambda_0+\mu_0}(m_g-m_0)\Big)\,U\,dr
   \ge 4\pi^2\int_0^\infty dr\int_0^\infty d\ell^2\int_{I_{r, \ell}^c} dw
   \,e^{\lambda_0}(1-\tilde{\eps}\,\sqrt{1+|v|^2})g. \]
In particular, if $g\in {\cal A}$ and ${\rm supp}\,g\subset {\rm supp}\,f_0$, then 
\begin{equation}\label{bayba} 
   \int_0^\infty\frac{d}{dr}\Big(e^{\lambda_0+\mu_0}(m_g-m_0)\Big)\,U\,dr\ge 0.
\end{equation} 
\end{cor}
{\bf Proof\,:} Define the functions
\[ a=\chi'(f_0)-1,\quad {\cal E}=\sqrt{1+|v|^2}=\sqrt{1+w^2+\ell^2/r^2},
   \quad b=\frac{a}{{\cal E}}=\frac{\chi'(f_0)-1}{\sqrt{1+w^2+\ell^2/r^2}}. \]   
We write 
\begin{eqnarray*} 
   & & \hspace{-3em}
   \int_0^\infty dr\int_0^\infty d\ell^2\int_{\R} dw
   \,e^{\lambda_0}\bigg[(\chi'(f_0)-1)(g-f_0)
   +e^{2\lambda_0}\,f_0\frac{d}{dw}\,[(\chi'(f_0)-1)w] 
   \,\frac{m_g-m_0}{r}\bigg]
   \\ & = & \int_0^\infty dr\int_0^\infty d\ell^2\int_{\R} dw
   \,e^{\lambda_0}\bigg[b\,{\cal E}(g-f_0)
   +e^{2\lambda_0}\,f_0\frac{d}{dw}\,(b {\cal E} w) 
   \,\frac{m_g-m_0}{r}\bigg]
   \\ & = & \int_0^\infty dr\int_0^\infty d\ell^2\int_{I_{r, \ell}} dw
   \,e^{\lambda_0}\bigg[b\,{\cal E}(g-f_0)
   +e^{2\lambda_0}\,f_0\frac{d}{dw}\,(b {\cal E} ¸w) 
   \,\frac{m_g-m_0}{r}\bigg]
   \\ & & +\,\int_0^\infty dr\int_0^\infty d\ell^2\int_{I_{r, \ell}^c} dw
   \,e^{\lambda_0} b\,{\cal E} g. 
\end{eqnarray*} 
On the interval $I_{r, \ell}$ we have $f_0(r, w, \ell)=(1-\tilde{\eps}(r)\sqrt{1+w^2+\ell^2/r^2})^k$, 
whence $\chi'(s)=s^{1/k}$ leads to $b(r, w, \ell)=-\tilde{\eps}(r)$, which is independent of $w$. 
On the other hand, on $I_{r, \ell}^c$ we have $b(r, w, \ell)=-1/{\cal E}$. 
This gives
\begin{eqnarray*} 
   & & \hspace{-3em}
   \int_0^\infty dr\int_0^\infty d\ell^2\int_{\R} dw
   \,e^{\lambda_0}\bigg[(\chi'(f_0)-1)(g-f_0)
   +e^{2\lambda_0}\,f_0\frac{d}{dw}\,[(\chi'(f_0)-1)w] 
   \,\frac{m_g-m_0}{r}\bigg]
   \\ & = & -\int_0^\infty dr\int_0^\infty d\ell^2\int_{I_{r, \ell}} dw
   \,e^{\lambda_0}\tilde{\eps}\bigg[{\cal E}(g-f_0)
   +e^{2\lambda_0}\,f_0\frac{d}{dw}\,({\cal E}w) 
   \,\frac{m_g-m_0}{r}\bigg]
   \\ & & -\,\int_0^\infty dr\int_0^\infty d\ell^2\int_{I_{r, \ell}^c} dw
   \,e^{\lambda_0}\,g
   \\ & = & -\int_0^\infty dr\int_0^\infty d\ell^2\int_{\R} dw
   \,e^{\lambda_0}\tilde{\eps}\,{\cal E}(g-f_0)
   +\int_0^\infty dr\int_0^\infty d\ell^2\int_{I_{r, \ell}^c} dw
   \,e^{\lambda_0}\tilde{\eps}\,{\cal E} g
   \\ & & -\int_0^\infty dr\int_0^\infty d\ell^2\int_{\R} dw
   \,e^{3\lambda_0}\,\tilde{\eps}\,f_0\frac{d}{dw}\,({\cal E}w)\,\frac{m_g-m_0}{r}
   \\ & & -\,\int_0^\infty dr\int_0^\infty d\ell^2\int_{I_{r, \ell}^c} dw
   \,e^{\lambda_0}\,g
   \\ & = & \int_0^\infty dr\,e^{\lambda_0+\mu_0}\,U\int_0^\infty d\ell^2\int_{\R} dw
   \,{\cal E}(g-f_0)
   \\ & & +\int_0^\infty dr\,e^{3\lambda_0+\mu_0}\,U\,\frac{m_g-m_0}{r}\int_0^\infty d\ell^2\int_{\R} dw
   \,f_0\Big({\cal E}+\frac{w^2}{{\cal E}}\Big)
   \\ & & -\int_0^\infty dr\int_0^\infty d\ell^2\int_{I_{r, \ell}^c} dw
   \,e^{\lambda_0}(1-\tilde{\eps}\,{\cal E})g
   \\ & = & \int_0^\infty dr\,e^{\lambda_0+\mu_0}\,U\,\frac{r^2}{\pi}\,(\rho_g-\rho_0)
   +\int_0^\infty dr\,e^{3\lambda_0+\mu_0}\,U\,\frac{m_g-m_0}{r}\,\frac{r^2}{\pi}(\rho_0+p_0)
   \\ & & -\int_0^\infty dr\int_0^\infty d\ell^2\int_{I_{r, \ell}^c} dw
   \,e^{\lambda_0}(1-\tilde{\eps}\,{\cal E})g. 
\end{eqnarray*}
Next note that, by definition of $\lambda_0$ and $\mu_0$, 
\begin{eqnarray*} 
   \lefteqn{\frac{d}{dr}\Big(e^{\lambda_0+\mu_0}(m_g-m_0)\Big)}
   \\ & = & 4\pi r^2\,e^{\lambda_0+\mu_0}(\rho_g-\rho_0)+e^{\lambda_0+\mu_0}(\lambda'_0+\mu'_0)(m_g-m_0)
   \\ & = & 4\pi r^2\,e^{\lambda_0+\mu_0}(\rho_g-\rho_0)
   +4\pi r^2\,e^{\lambda_0+\mu_0}e^{2\lambda_0}\,\frac{m_g-m_0}{r}\,(\rho_0+p_0), 
\end{eqnarray*} 
which yields the claim, if we use Corollary \ref{reform}. 
{\hfill$\Box$}\bigskip

Clearly, if it can be shown that $U(r)={\rm const.}$, then the minimizer $f_0$ is a static solution having the form
\begin{equation}\label{fastform2} 
   f_0(r, w, \ell)
   =\Big(1-Ce^{\mu_0(r)}\sqrt{1+w^2+\ell^2/r^2}\Big)_+^k. 
\end{equation}
The following lemma implies that if the induced density stays away from the boundary 
in the set $\mathcal{A}$ and if $f_0>0$ on a strip, then $f_0$ takes the form (\ref{fastform2}) in this domain. 

\begin{lemma}\label{lemmaUconst} 
Suppose that $\rho_0(r)\le\sigma_0-\delta_0$ for some $\delta_0>0$. If 
\[ f_0(r, w, \ell)\ge\eps_0\quad\mbox{a.e.~in}\quad [r_1, r_2]\times [w_1, w_2]\times [\ell_1, \ell_2] \] 
for some $\eps_0>0$, $0<r_1<r_2<R_0$, $w_1\le w_2$ and $0\le\ell_1\le\ell_2$, then 
\[ U(r)={\rm const.}\quad\mbox{a.e.~in}\quad [r_1, r_2]. \] 
\end{lemma} 
{\bf Proof\,:} For $A>0$ consider 
\[ g=f_0+\frac{A}{\sqrt{1+|v|^2}}\,{\bf 1}_{\{|v|\le 1\}}\,\frac{\eta'(r)}{4\pi r^2}
   \,{\bf 1}_{[w_1, w_2]}(w)\,{\bf 1}_{[\ell_1, \ell_2]}(\ell), \] 
where $\eta\in C^1([0, \infty[)$ is such that $\eta(r)=0$ outside of $[r_1, r_2]$. 
Since $r_1>0$ we have ${\|\frac{\eta'(r)}{4\pi r^2}\|}_{L^\infty_r}\le C_0$. 
If $(r, w, \ell)\in Q=[r_1, r_2]\times [w_1, w_2]\times [\ell_1, \ell_2]$, 
then $f_0(r, w, \ell)\ge\eps_0$. Thus if $|A|\le\frac{\eps_0}{2C_0}$, 
then $g(r, w, \ell)\ge 0$ in $Q$. On the other hand, for $(r, w, \ell)\not\in Q$ 
we get $g(r, w, \ell)=f_0(r, w, \ell)\ge 0$. Integrating against $\sqrt{1+|v|^2}$ in $v$, 
we obtain 
\[ \rho_g(r)=\rho_0(r)+AC_1\,\frac{\eta'(r)}{4\pi r^2} \] 
for $C_1=\int {\bf 1}_{\{|v|\le 1\}}\,{\bf 1}_{[w_1, w_2]}(w)\,{\bf 1}_{[\ell_1, \ell_2]}(\ell)\,dv$. 
Hence if $|A|\le\frac{\delta_0}{C_0 C_1}$, then $\rho_g(r)\le\sigma_0$ for all $r$.  
For the mass functions it follows 
\[ m_g(r)=m_0(r)+AC_1\,\int_0^r\eta'(s)\,ds=m_0(r)+AC_1\eta(r), \] 
and in particular $\eta(r)=0$ for $r\ge r_2$ implies that $g$ has total mass $M$. 
In summary, we have shown that $g\in {\cal A}$. 
If $f_0(r, w, \ell)=0$, then $(r, w, \ell)\not\in Q$ 
and $g(r, w, \ell)=f_0(r, w, \ell)=0$. Thus (\ref{bayba}) is applicable and yields 
\[ \int_0^\infty\frac{d}{dr}(e^{\lambda_0+\mu_0}\eta)\,U\,dr\ge 0. \] 
Since $\eta$ can be replaced by $-\eta$, it follows that 
\begin{equation}\label{hauhe} 
   \int_0^\infty\frac{d}{dr}(e^{\lambda_0+\mu_0}\eta)\,U\,dr=0
\end{equation} 
for all $C^1$-function $\eta$ that vanish outside of $[r_1, r_2]$. 
Let $\varphi\in C_0^\infty(]r_1, r_2[)$ be a test function. By extending $\varphi$ 
by zero outside of $[r_1, r_2]$, we consider $\varphi$ to be defined on $[0, \infty[$. 
Then $\eta=e^{-(\lambda_0+\mu_0)}\varphi$ is a $C^1$-function that vanishes outside of $[r_1, r_2]$. 
From (\ref{hauhe}) we therefore deduce 
\[ \int_{r_1}^{r_2}\varphi'\,U\,dr=0 \] 
for all $\varphi\in C_0^\infty(]r_1, r_2[)$, so that $U(r)={\rm const.}$ a.e.~in $[r_1, r_2]$. 
{\hfill$\Box$}\bigskip

Let 
\[ R_0=\inf\{R>0: f_0(r, w, \ell)=0\,\,\mbox{for a.e.}
   \,\,r\in [R, \infty[, w\in\R, \ell\in [0, \infty[\}. \] 

Next we show a further support property of $f_0$. 

\begin{lemma}\label{f0-no-gaps}
Suppose that $\rho_0(r)\le\sigma_0-\delta_0$ for some $\delta_0>0$. 
Then there is no interval $]a, b[$ such that $0<a<b<R_0$ and $f_0(r, w, \ell)=0$ 
for all a.e.~$(r, w, \ell)\in ]a, b[\times\R\times [0, \infty[$. 
\end{lemma}
{\bf Proof\,:} Assume, on the contrary, that there is such an interval $]a, b[$. 
By moving $b$ closer to $a$ if necessary, 
we may suppose that $\delta=b-a\le\min\{\frac{1}{2}R_0, ((\frac{\sigma_0}{\sigma_0-\delta_0})^{1/2}-1)a\}$. 
Define $f_0^-:=f_0\,\textbf{1}_{[0, a]}(|x|)$ and $f_0^+:=f_0\,\textbf{1}_{[b, R_0]}(|x|)$ 
to obtain $f_0=f_0^- +f_0^+$. For $r\in [a, \infty[$ define $f^+(r, w, \ell)=f_0^+(r+\delta, w, \ell)$, 
i.e., $f^+=f_0 \textbf{1}_{[a, R_0-\delta]}(|x|)$, and furthermore $f=f_0^- + f^+$. We also put 
\[ F_0(r)=\int_0^\infty d\ell^2\int_{\R} dw\,\sqrt{1+w^2+\ell^2/r^2}\,f_0(r, w, \ell) \]
and
\[ G_0(r)=\int_0^\infty d\ell^2\int_{\R} dw\,\Big(\chi(f_0(r, w, \ell))-f_0(r, w, \ell)\Big); \]
the functions $F$ and $G$ are defined analogously by replacing $f_0$ by $f$. Note that 
\begin{equation}\label{1919} 
   \rho_0(r)=\frac{\pi}{r^2}\,F_0(r),\quad\rho(r)=\frac{\pi}{r^2}\,F(r),
\end{equation}  
if $\rho$ denotes the density as induced by $f$. If $r\in [0, a]$, 
then $f_0^+(r, \cdot, \cdot)=f^+(r, \cdot, \cdot)=0$ a.e.~and hence $F_0(r)=F(r)$ as well as $G_0(r)=G(r)$. 
If $r\in [a, \infty[$, then $f_0^-(r, \cdot, \cdot)=0$ a.e., so that $F_0(r+\delta)=F(r)$ 
and $G_0(r+\delta)=G(r)$ by definition. For the mass functions, if $r\in [0, a]$, 
then $m_0(r)=m(r)$ by (\ref{1919}). If $r\in [a, \infty[$, then 
\begin{eqnarray*}
   m(r) & = & m(a)+4\pi^2\int_a^r F(s)\,ds
   =m_0(a)+4\pi^2\int_a^r F_0(s+\delta)\,ds
   \\ & = & m_0(b)+4\pi^2\int_{b}^{r+\delta} F_0(s)\,ds=m_0(r+\delta).
\end{eqnarray*}
Here we used that $m_0(a)=m_0(b)$ since $f_0=0$ on $]a, b[$. 
In particular, the induced total masses of $f_0$ and $f$ are the same. 
Moreover, 
\[ \rho(r)=\frac{\pi}{r^2}\,F(r)
   =\left\{\begin{array}{c@{\quad:\quad}l}
   \frac{\pi}{r^2}\,F_0(r) & r\in [0, a] 
   \\[1ex] \frac{\pi}{r^2}\,F_0(r+\delta) & r\in [a, R_0-\delta] 
   \\[1ex] 0 & r\in [R_0-\delta, \infty[\end{array}\right.
   =\left\{\begin{array}{c@{\quad:\quad}l}
   \rho_0(r) & r\in [0, a] 
   \\[1ex] \frac{(r+\delta)^2}{r^2}\,\rho_0(r+\delta) & r\in [a, R_0-\delta] 
   \\[1ex] 0 & r\in [R_0-\delta, \infty[\end{array}\right. . \] 
If $r\ge a$, then $\frac{(r+\delta)^2}{r^2}\rho_0(r+\delta)\le\sigma_0$ 
by the choice of $\delta$ and the assumption on $\rho_0$. 
Therefore we obtain $\rho(r)\le\sigma_0$ for all $r\in [0, \infty[$, 
and in turn we have shown that $f\in {\cal A}$. Lastly, we observe that 
\begin{eqnarray*}
   {\cal D}(f) & = & 4\pi^2\int_0^{R_0-\delta} \frac{G(r)}{\sqrt{1-\frac{2m(r)}{r}}}\,dr
   \\ & = & 4\pi^2 \int_0^a \frac{G_0(r)}{\sqrt{1-\frac{2m_0(r)}{r}}}\,dr
   +4\pi^2 \int_a^{R_0-\delta}\frac{G(r)}{\sqrt{1-\frac{2m(r)}{r}}}\,dr=:T_1+T_2. 
\end{eqnarray*}
To estimate $T_2$, we recall from Theorem \ref{f0exists} that $f_0\le 1$, 
which implies that $G_0$ is non-positive. It follows that 
\begin{eqnarray}\label{wini21}
   T_2 & = & 4\pi^2 \int_a^{R_0-\delta}\frac{G_0(r+\delta)}{\sqrt{1-\frac{2m_0(r+\delta)}{r}}}\,dr
   \leq 4\pi^2 \int_a^{R_0-\delta}\frac{G_0(r+\delta)}{\sqrt{1-\frac{2m_0(r+\delta)}{r+\delta}}}\,dr
   \nonumber   
   \\ & = & 4\pi^2 \int_b^{R_0}\frac{G_0(r)}{\sqrt{1-\frac{2m_0(r)}{r}}}\,dr.
\end{eqnarray}
Actually if we had equality in (\ref{wini21}), then $G_0(r+\delta)m_0(r+\delta)=0$ 
a.e.~in $[a, R_0-\delta]$, which means that $G_0(r)m_0(r)=0$ a.e.~in $[b, R_0]$. 
If there is a $r_1\in [b, R_0]$ such that $m_0(r_1)=0$, then there is also a last $r_0\in [b, R_0]$ 
with this property; recall that $m_0$ is continuous. From $m_0(r_0)=0$ it then follows that 
$f_0=0$ a.e.~in $[0, r_0]\times\R\times [0, \infty[$. Now $r_0=R_0$ is impossible, 
since this would yield $f_0=0$ a.e., contradicting $I<0$ from Lemma \ref{Ikle}. 
Therefore we have $r_0<R_0$ and thus $G_0(r)=0$ for a.e.~$r\in [r_0, R_0]$, 
which by virtue of $\chi(f_0)-f_0\le 0$ in turn yields $f_0=0$ a.e.~in $[r_0, R_0]$, 
and yet a further contradiction, this time to the definition of $R_0$. 
As a consequence, $m_0(r)\neq 0$ for $r\in [b, R_0]$, 
so that $G_0(r)=0$ for a.e.~$r\in [b, R_0]$ and hence $f_0=0$ a.e.~in $[b, R_0]\times\R\times [0, \infty[$, 
which is again impossible by the definition of $R_0$. 
Thus indeed $T_2<4\pi^2 \int_b^{R_0} (1-\frac{2m_0(r)}{r})^{-1/2}\,G_0(r)\,dr$ in (\ref{wini21}). 
Hence we obtain
\[ {\cal D}(f)<4\pi^2 \int_0^a\frac{G_0(r)}{\sqrt{1-\frac{2m_0(r)}{r}}}\,dr 
   +4\pi^2\int_b^{R_0}\frac{G_0(r)}{\sqrt{1-\frac{2m_0(r)}{r}}}\,dr={\cal D}(f_0). \] 
Since $f_0$ is a minimizer, this is a contradiction. 
{\hfill$\Box$}\bigskip

%%%%%%%%%%%%%%%%%%%%%%%%%%%%%%%%%%%%%%%%%%%%%%%%%%%%%%%%%%%%%%%%%%%%%%%%%%%%%%%    

%\section{Small mass}
%\label{smama} 
%\setcounter{equation}{0}

In the following lemma we show that there are minimizers that are not constant $\rho_0(r)=\sigma_0$ 
on the entire interval $[0, M/\beta]$, if the mass $M$ is small. 
This indicates that in this case one ultimately might be able 
to prove that $\rho_0(r)$ stays below $\sigma_0$ on the whole interval, 
which is in accordance with what we have seen numerically. 

\begin{lemma}\label{entire}
Let $\beta$ be sufficiently small so that $e^{\lambda}<21/20$. Let $M=\sqrt{3\beta^3/4\pi}$ 
so that $\sigma_0=1$, cf.~(\ref{siggi}). If $\rho_0=\sigma_0=1$ on $[0,r_\ast]$, then $r_\ast<M/\beta$. 
\end{lemma}
{\bf Proof\,:} 
Assume that $\rho_0=\sigma_0=1$ on $[0, M/\beta]$. 
Let $K=12/25$. Since $f_0\leq 1$ we have
\[ \int_{|v|\leq K}\sqrt{1+|v|^2}f_0\, dv\leq \int_{|v|\leq K}\sqrt{1+|v|^2}\, dv\leq \frac{\sigma_0}{2}=\frac12, \]
which implies that
\[ M_{<}=4\pi\int_0^{M/\beta}r^2\int_{|v|\leq K}\sqrt{1+|v|^2}f_0\,dv\,dr\leq\frac{M}{2}. \]
Clearly, we have that $M_{>}=M-M_{<}$. From (\ref{rbg2}) we get
\[ |\mathcal{D}_{>}|\leq \frac{21}{20}\,\frac{M_{>}}{\sqrt{1+K^2}}
   =\frac{21}{20}\,\frac{(M-M_{<})}{\sqrt{1+K^2}}. \]
In addition, from (\ref{rbg1}) we have
\[ |\mathcal{D}_{<}|\leq\frac{21}{20}\,M_<\leq\frac{21}{40}\,M. \]
Hence we obtain 
\[ |\mathcal{D}|\leq \frac{21}{20}\,M_<+\frac{21}{20}\,\frac{(M-M_{<})}{\sqrt{1+K^2}}
   \leq\frac{21}{20}\,\frac{M}{\sqrt{1+K^2}}+\frac{21}{40}\,M\Big(1-\frac{1}{\sqrt{1+K^2}}\Big)<0.9983\,M<M, \]
which contradicts (CBEC).
{\hfill$\Box$}\bigskip

\begin{remark}
In the argument above we used that $f_0\leq 1$. In the case of small mass, one can deduce that $f_0\ll 1$, at least on a large set of its support. This fact can be used to improve the result, so that $\rho< \sigma_0$ on a large part of the interval, but since we anyway have not been able to prove this on the entire interval we have settled for the result above. 
\end{remark}

%%%%%%%%%%%%%%%%%%%%%%%%%%%%%%%%%%%%%%%%%%%%%%%%%%%%%%%%%%%%%%%%%%%%%%%%%%%%%%%    

\section{On the Casimir-binding-energy condition (CBEC)}
\label{CBECdisc_sect} 
\setcounter{equation}{0}

In this section we discuss the relevance of the Casimir-binding-energy condition (CBEC) for static solutions, we find strong evidence that static solutions that are not too relativistic will satisfy this condition. 
We also comment on the stability conjecture proposed in the introduction. 

Let us start by deriving a couple of general results for static solutions. 

\begin{lemma}\label{Hf0}
Suppose that $f_0$ is a static solution of the form (\ref{fastform2}), i.e., 
$f_0=(1-e^{\mu_0(r)-\mu_0(R_0)}{\cal E})_+^k$ for ${\cal E}=\sqrt{1+|v|^2}=\sqrt{1+w^2+\ell^2/r^2}$. Then 
\[ {\cal D}(f_0)=-4\pi e^{-\mu_0(R_0)}\int_0^{R_0} r^2 (p_0+\rho_0)\,e^{\mu_0+\lambda_0}\,dr. \]
\end{lemma}
{\bf Proof\,:} From the proof of Corollary \ref{reform} it follows that
\begin{eqnarray*}
  \int_{\R^3} (\chi (f_0)-f_0)\,dv & = & \frac{\pi}{r^2}\int_0^\infty d\ell^2\int_{\R} dw\,
  \frac{d}{dw}\,[(\chi'(f_0)-1)w]\,f_0 
  \\ & = &  -\frac{\pi}{r^2}\,e^{\mu_0(r)-\mu_0(R_0)}\int_0^\infty d\ell^2\int_{\R} dw\,
  \frac{d}{dw}\,[{\cal E} w]\,f_0 
  \\ & = &  -\frac{\pi}{r^2}\,e^{\mu_0(r)-\mu_0(R_0)}\int_0^\infty d\ell^2\int_{\R} dw\,
  \Big(\frac{w^2}{\cal E}+{\cal E}\Big)\,f_0 
  \\ & = & -(p_0+\rho_0)\,e^{\mu_0(r)-\mu_0(R_0)}.
\end{eqnarray*}
The claim follows by inserting this expression into the formula for ${\cal D}$. 
{\hfill$\Box$}\bigskip

\begin{lemma}\label{Hf0bounds}
Suppose that $f_0=(1-e^{\mu_0(r)-\mu_0(R_0)}{\cal E})_+^k$. Then
\[ \frac12\,e^{-\mu_0(R_0)}M\leq|{\cal D}(f_0)|\leq e^{-\mu_0(R_0)}M. \]
\end{lemma}
{\bf Proof\,:} From \cite[equ.~(10)]{A1} it follows that 
\begin{equation}\label{M10}
4\pi \int_0^{R_0} r^2 (\rho_0+p_0+q_0)e^{\mu_0+\lambda_0}\, dr=M. 
\end{equation}
Here $q_0$ denotes twice the tangential pressure. 
From this relation and Lemma \ref{Hf0} we can derive bounds for $|{\cal D}(f_0)|$. We immediately have
\[ |{\cal D}(f_0)|\leq e^{-\mu_0(R_0)}M, \]
and furthermore, since $\rho_0\geq p_0+q_0$,
\begin{eqnarray*}
   |{\cal D}(f_0)|& \geq & 4\pi e^{-\mu_0(R_0)}\int_0^{R_0} r^2 \rho_0\,e^{\mu_0+\lambda_0}\,dr
   \\ & \geq & 2\pi e^{-\mu_0(R_0)}\int_0^{R_0} r^2 (\rho_0+p_0+q_0)\,e^{\mu_0+\lambda_0}\,dr
   \\ & = & \frac12\,e^{-\mu_0(R_0)}M, 
\end{eqnarray*}
which completes the argument. 
{\hfill$\Box$}\bigskip

We note that in our case the pressure is isotropic and we have $q_0=2p_0$, so that 
\[ 4\pi \int_0^{R_0} r^2 (\rho_0+3 p_0)\,e^{\mu_0+\lambda_0}\,dr=M.  \]
From Lemma {\ref{Hf0bounds} it first seems that highly relativistic solutions for which 
\[ e^{-\mu_0(R_0)}=\frac{1}{1-\frac{2M}{R_0}} \] 
is large should be good candidates to satisfy the (CBEC). Note that $2M/R_0<8/9$ 
for any static solution of the Einstein-Vlasov system, cf.~\cite{A3}, 
whereas for isotropic solutions an analytic bound is not known, but numerical simulations indicate that $2M/R_0<1/2$ 
in that case. We now argue why it is misleading that highly relativistic solutions are the best candidates. 
For, let us consider highly relativistic isotropic solutions where the extreme case $2M/R_0=1/2$ is attained. 
Then Lemma {\ref{Hf0bounds} gives the lower bound
\[ |{\cal D}(f_0)|\geq \frac12\,e^{-\mu_0(R_0)}M=\frac{M}{\sqrt{2}}. \]
On the other hand, for static solutions that are not very relativistic, 
but close to being Newtonian, we have $p_0\ll\rho_0$, so that
\begin{eqnarray*}
   |{\cal D}(f_0)| & = & 4\pi e^{-\mu_0(R_0)}\int_0^{R_0} r^2 (p_0+\rho_0)\,e^{\mu_0+\lambda_0}\,dr\\
   && \simeq e^{-\mu_0(R_0)}4\pi \int_0^{R_0} r^2 (\rho_0+3 p_0)\,e^{\mu_0+\lambda_0}\,dr=e^{-\mu_0(R_0)}M>M.
\end{eqnarray*}
Hence a static solution that is not too relativistic is a good candidate for satisfying the (CBEC). 
Now, it is straightforward to numerically check (CBEC) on particular static solutions, and indeed 
we find that solutions that are not too relativistic do verify the condition. 

In the introduction we conjecture that isotropic static solutions which satisfy the (CBEC) are stable with respect to mass-preserving perturbations. 
%We have found numerical support for this statement although the numerical study we have done so far is very limited. 
An interesting observation is that when $k\ll 1$, the Casimir-binding energy $E_{Cb}$ and the binding energy $E_b$ almost coincide. If the conjecture is true it indicates that when $k=0$, the sign of $E_b$ alone determines whether a static solution is stable or unstable. The first author investigated, together with Gerhard Rein, stability of static shells in \cite{AR}. (Note that the conjecture we pose concerns isotropic solutions rather than shells, which are anisotropic, but a comparison is nevertheless relevant.) In \cite{AR} the perturbations had the form $f=Af_0$, where $A<1$ or $A>1$, which we here call amplitude-perturbations. From the results in Table 1 in \cite{AR} it follows that the sign of $E_b$ does not determine stability in the case of amplitude-perturbations; if $A>1$ there are several cases when the static solution is unstable even if $E_b>0$. (The binding energy in \cite{AR} is in fact the fractional binding energy but the sign is the same as of $E_b$.) Now, the case $A>1$ does not preserve the mass and such a perturbation is outside the scope of our conjecture. In a very limited study we have investigated the case $k=0$ for shells and we \textit{do} find numerical support that for mass-preserving perturbations the sign of $E_b$ alone determines the stability properties. An interesting aspect of this is that in the numerical study \cite{Getal}, another type of perturbations were considered which are called dynamically accessible perturbations. The results in \cite{Getal} indicate that the stability properties for dynamically accessible perturbations and amplitude-perturbations are similar in the following sense: if an amplitude-perturbation with $A>1$ is unstable then a dynamically accessible perturbation is also unstable if the perturbation is such that the particles of the perturbation initially go inward, i.e. the steady state is pushed towards collapse rather than towards dispersion. 
%Hence, they conclude that the decisive property is if the steady state is pushed towards collapse or towards dispersion.
In our situation with mass-preserving perturbations we have therefore tried perturbations of both types, i.e. perturbations that push towards collapse and towards dispersion. In both cases we find numerical support that a shell with $k=0$, and such that $E_b>0$, is stable in accordance with our conjecture. This indicates an essential difference between mass-preserving perturbations and dynamically accessible perturbations. A more extensive numerical study will be carried out in future work to fully understand this intriguing issue. 

%We also observe that these solutions 
%are located before the first local maximum of the binding energy along a sequence of static solutions,  
%which supports the belief that they are stable in view of the weak binding energy hypothesis. 

%%%%%%%%%%%%%%%%%%%%%%%%%%%%%%%%%%%%%%%%%%%%%%%%%%%%%%%%%%%%%%%%%%%%%%%%%%%%%%%%%%%%%%%%

\section{Some technical results} 
\label{tech_sect} 
\setcounter{equation}{0}

\subsection{A convexity property}

The following observation is a crucial ingredient of the minimization argument, 
and this fact seems not to have been noticed so far.

\begin{lemma}\label{phikonv} Let $k\in ]0, 2]$. Then the functional
\[ f\mapsto\Psi(f)=\int_{\R^3}\int_{\R^3} e^{\lambda}\,\chi(f)\,dx\,dv \]
is convex.
\end{lemma}
{\bf Proof\,:} Let $\varphi\in C_0^\infty(\R^3\times\R^3)$. To begin with,
\[ e^{\lambda_f}=\frac{1}{\sqrt{1-\frac{2m}{r}}}. \]
Since $m=m_f$ is linear in $f$, we obtain
\[ D(e^{\lambda_f})(\varphi)=\frac{1}{r}\,\Big(1-\frac{2m_f}{r}\Big)^{-3/2} m_\varphi
   =\frac{1}{r}\,e^{3\lambda_f}\,m_\varphi \]
for the derivative. Thus
\[ D\Psi(f)(\varphi)=\int_{\R^3}\int_{\R^3}\bigg(\frac{1}{r}\,e^{3\lambda_f}
   \,m_\varphi\,\chi(f)+e^{\lambda_f}\,\chi'(f)\varphi\bigg)\,dx\,dv. \]
In addition, $D(e^{3\lambda_f})(\varphi)=\frac{3}{r}\,e^{5\lambda_f}\,m_\varphi$. 
Hence, for the second derivative,
\begin{equation}\label{D2Phi}
   D^2\Psi(f)(\varphi, \varphi)=\int_{\R^3}\int_{\R^3}\bigg[\frac{3}{r^2}\,
   e^{5\lambda_f}\,m_\varphi^2\,\chi(f)
   +\frac{2}{r}\,e^{3\lambda_f}\,m_\varphi\,\chi'(f)\varphi
   +e^{\lambda_f}\,\chi''(f)\varphi^2\bigg]\,dx\,dv.
\end{equation}
Using $\chi(f)=\frac{k}{k+1}\,f^{1+1/k}$, $\chi'(f)=f^{1/k}$ 
and $\chi''(f)=\frac{1}{k}\,f^{1/k-1}$,
the relation (\ref{D2Phi}) can be rewritten as
\begin{eqnarray*}
   \lefteqn{D^2\Psi(f)(\varphi, \varphi)}
   \\ & = & \frac{k}{k+1}\int_{\R^3}\int_{\R^3} e^{\lambda_f} f^{1/k-1}
   \bigg[\frac{3}{r^2}\,e^{4\lambda_f}\,m_\varphi^2\,f^2
   +\frac{2}{r}\,e^{2\lambda_f}\,m_\varphi\,\Big(1+\frac{1}{k}\Big)f\varphi
   +\frac{1}{k}\Big(1+\frac{1}{k}\Big)\varphi^2\bigg]\,dx\,dv
   \\ & = & \frac{k}{k+1}\int_{\R^3}\int_{\R^3} e^{\lambda_f} f^{1/k-1}
   \bigg[\bigg(\frac{1}{\sqrt{3}}\Big(1+\frac{1}{k}\Big)\varphi
   +\sqrt{3}\,\frac{1}{r}\,e^{2\lambda_f} m_\varphi f\bigg)^2
   +\frac{(2-k)}{3k}\Big(1+\frac{1}{k}\Big)\varphi^2\bigg]\,dx\,dv,
\end{eqnarray*}
which is positive, owing to $2-k\ge 0$. Thus $\Psi$ is convex.
{\hfill$\Box$}\bigskip

\subsection{A lemma on $2m/r$}

\begin{lemma}\label{2mrlem} 
Let $M=\int_{\R^3}\int_{\R^3}\sqrt{1+|v|^2}\,f\,dx\,dv$ be the mass and 
\[ \sigma_M=\frac{3}{32\pi M^2}. \] 
If $\rho(r)\le\sigma\le\sigma_M$ for all $r$, then 
\[ \frac{2m(r)}{r}\le\bigg(\frac{\sigma}{\sigma_M}\bigg)^{1/3}
   \quad\mbox{for all}\quad r. \] 
In particular, if $\sigma<\sigma_M$, then 
\[ \frac{2m(r)}{r}\le\bigg(\frac{\sigma}{\sigma_M}\bigg)^{1/3}<1
   \quad\mbox{for all}\quad r. \]
Furthermore, if $\rho(r)\le\sigma_M$ for all $r$ and
\begin{equation}\label{betpa} 
   \frac{2m(r)}{r}=1\quad\mbox{for some}\quad r=r_\ast,
\end{equation}  
then $r_\ast=2M$. 
\end{lemma}    
{\bf Proof\,:} To establish the first claim, 
let $r_1=(\frac{3M}{4\pi\sigma})^{1/3}$. From the assumption we have 
\[ m(r)=4\pi\int_0^r \eta^2\,\rho(\eta)\,d\eta
   \le\frac{4\pi}{3}\,\sigma\,r^3 \]
for all $r$. Hence, if $r\le r_1$, then we get 
\begin{equation}\label{spri1} 
   \frac{2m(r)}{r}\le\frac{8\pi}{3}\,\sigma\,r^2
   \le\frac{8\pi}{3}\,\sigma\,r_1^2
   =\frac{8\pi}{3}\,\sigma\,\bigg(\frac{3M}{4\pi\sigma}\bigg)^{2/3}
   =\bigg(\frac{\sigma}{\sigma_M}\bigg)^{1/3}.
\end{equation}  
On the other hand, if $r\ge r_1$, then also 
\begin{equation}\label{spri2} 
   \frac{2m(r)}{r}\le\frac{2M}{r}\le\frac{2M}{r_1}
   =2M\bigg(\frac{4\pi\sigma}{3M}\bigg)^{1/3}
   =\bigg(\frac{\sigma}{\sigma_M}\bigg)^{1/3}.
\end{equation} 
For the second assertion, suppose that $\rho(r)\le\sigma_M$ for all $r$ 
and that (\ref{betpa}) holds. We are going to take $\sigma=\sigma_M$ 
in the first part of the argument. If we had $r_\ast<r_1$, 
then repeating (\ref{spri1}) with the strict inequality 
would imply 
\[ 1=\frac{2m(r_\ast)}{r_\ast}<1, \]  
which is impossible. Accordingly, we must have $r_\ast\ge r_1$. But then 
(\ref{spri2}) leads to 
\[ 1=\frac{2m(r_\ast)}{r_\ast}\le\frac{2M}{r_\ast}\le 1, \] 
which gives $r_\ast=2M$ as desired. 
{\hfill$\Box$}\bigskip

%%%%%%%%%%%%%%%%%%%%%%%%%%%%%%%%%%%%%%%%%%%%%%%%%%%%%%%%%%%%%%%%%%%%%%%%%%%%%%%%%%%%%%%%%

\end{document}